\documentclass[10pt,a4paper]{article}
\usepackage{amscd}

\usepackage{mathrsfs}
\usepackage{amsfonts}
\usepackage{amssymb}
\usepackage{pifont}
\usepackage{amsmath}
\usepackage{mathrsfs}
\usepackage{enumerate}
\usepackage{amsmath}
\usepackage{mathrsfs}
\usepackage{amsfonts}
\usepackage{amsmath}
\usepackage{amscd}
\usepackage[all]{xy}

\usepackage{amsthm}
\usepackage{amsmath}
\usepackage{mathrsfs}
\usepackage{amsfonts}
\usepackage{mathrsfs}
\usepackage{amsfonts}
\usepackage{mathrsfs}
\usepackage{bbding}
\usepackage{amsfonts}
\usepackage{tikz}
\usetikzlibrary{%
  calc,
  fit,
  shapes,
  backgrounds
}

\renewcommand\appendix{\par
  \setcounter{section}{0}
  \setcounter{subsection}{0}
  \setcounter{figure}{0}
  \setcounter{table}{0}
  \renewcommand\thesection{Appendix \Alph{section}}
  \renewcommand\thefigure{\Alph{section}\arabic{figure}}
  \renewcommand\thetable{\Alph{section}\arabic{table}}
}

\usepackage{txfonts}
\usepackage{latexsym}
\usepackage{amssymb}
\usepackage{amstext}
\usepackage{fancyhdr}
\usepackage{graphicx,picinpar,epsf}
\usepackage{lastpage}
\usepackage{hyperref}
\usepackage{verbatim}
\newtheorem{theorem}{Theorem}[section]
\newtheorem{definition}[theorem]{Definition}
\newtheorem{proposition}[theorem]{Proposition}

\newtheorem{corollary}[theorem]{Corollary}
\newtheorem{lemma}[theorem]{Lemma}
\newtheorem{remark}[theorem]{Remark}
\newtheorem{example}[theorem]{Example}

\title{The MMP for deformations of Hilb$^n$ \textbf P$^2$}
\author{Chunyi Li\\
University of Illinois\\
\texttt{li118@illinois.edu}
 \and Xiaolei Zhao\\
 University of Michigan\\
 \texttt{xlzhao@umich.edu}}
\date{\today}

\begin{document}
\maketitle
\begin{abstract}
We study the birational geometry of deformations of Hilbert schemes
of points on \textbf P$^2$. On one hand, we complete the picture in
\cite{ABCH} by giving an explicit correspondence between the stable
base locus walls on the Neron-Severi space and the actual walls on
the Bridgeland stability space. On the other hand, we show that the
birational geometry of a deformed Hilb P$^2$ is different from that
of Hilb P$^2$.
\end{abstract}

\section*{Introduction}
The Hilbert scheme of points on an algebraic variety is the moduli
space that parameterizes all of the $0$-dimensional subschemes with
length $n$ on the variety, where $n$ is a given positive integer. In
the case when the variety is a curve, the Hilbert scheme is just the
symmetric product of the curve itself. While the variety has
dimension greater than $2$, the Hilbert scheme has very bad
singularities. In the surface case, the Hilbert scheme of points is
smooth and connected which becomes a nice geometric object for
study. As a moduli space of sheaves on the surface, the Hilbert
scheme parameterizes the ideal sheaves with trivial first Chern
class and a given second Chern class. The goal of this paper is to
understand the birational geometry of the Hilbert scheme and the
deformed Hilbert scheme of points on the projective plane. This sits
into a huge program which studies the birational geometry of moduli
spaces. In the case of the moduli space of curves, this is known as
the Hassett-Keel program, and much research has been done. In the
case of the moduli space of sheaves, much recent progress is made
via Bridgeland stability conditions.\\

The notation of the stability condition on a triangulated category T
has been introduced by Bridgeland in \cite{Br}. It is given by
abstracting the usual properties of the $\mu$-stability for sheaves
on projective varieties. The central charge, which substitutes the
slope $\mu$, is a group homomorphism from the numerical Grothendick
group to $\mathbb C$, and satisfies some extra conditions including
the Harder-Narasimhan filtrations. The Bridgeland stability
conditions form a natural topological space Stab(T), which becomes a
manifold of dimension not exceeding that of
the numerical Grothendieck group.\\
Consider the case when T is a bounded derived category of coherent
sheaves on a smooth surface X. Given a numerical equivalence class
\textbf v, and a stability condition $\sigma$ on $D^b(Coh(X))$, we
have the moduli space $\mathfrak M_{\sigma}(\textbf v)$ of
$\sigma$-stable complexes of numerical type \textbf v. Two natural
questions on $\mathfrak M_{\sigma}(\textbf v)$ are:
\begin{enumerate}
\item When is $\mathfrak M_{\sigma}(\textbf v)$ a good geometric object to
study?
\item When $\sigma$ changes in Stab$(D^b(X))$, what is the behavior
of $\mathfrak M_{\sigma}(\textbf v)$?
\end{enumerate}

For a general smooth surface, the known answers to both questions
are mostly either vague or philosophical. For the first question,
there are a few ways to determine when $\mathfrak M_{\sigma}(v)$ is
projective. For the second question, ideally, the stability space
has a well-behaved chamber structure. In each chamber, $\mathfrak
M_\sigma(v)\simeq \mathfrak M_{\sigma '}(v)$. Among different
chambers, there is birational map $\mathfrak M_\sigma(v)
\dashrightarrow \mathfrak M_{\sigma '}(v)$. Yet the ideal picture is
far from being accomplished. It is only set up or
partially/conjecturally set up when $X$ is a K3 surface, the
projective plane, a high degree del Pezzo surface, a Hirzebruch
surface or an abelian surface.\\

Let $\mathfrak M_\sigma(n)$ be the moduli space of complexes on
\textbf P$^2$ with numerical type $(r,c_1,\chi)$ $=$ $(1,0,1-n)$,
i.e., the numerical type of Hilbert schemes. In \cite{ABCH}, the
authors studied the two questions in this case. They describe a wall
and chamber structure on Stab$(D^b(\textbf P^2))$ for the invariant
$(r,c_1,\chi)$ $=$ $(1,0,1-n)$. On a particular upper half plane
slice of Stab$(D^b (\textbf P^2))$, the walls are a sequence of
nested semicircles in each quadrant, plus the vertical axis. For
certain $\sigma$ in the second quadrant, they show that $\mathfrak
M_{\sigma}(n)$ $\simeq$ Hilb$^n\textbf P^2$. By choosing certain
representative stability condition in each chamber, the authors
prove that there are finitely many chambers for which $\mathfrak
M_{\sigma}(n)$ is non-empty, and, in this case, projective.
Moreover, for small values of $n$, the author also write down an
explicit correspondence between the chamber walls of the stability
space and the base locus decomposition walls of the effective cone
(in the sense of MMP). For general value of n, this explicit formula
of the correspondence remains conjectural. One difficulty is to get
a better answer to the first question above, in other words, to
control the behavior of $\mathfrak M_{\sigma}(n)$, especially the
smoothness and irreducibility. In this paper, we solve these two
questions in this case of \textbf P$^2$ with numerical type
$(r,c_1,\chi)$ $=$ $(1,0,1-n)$.

\begin{theorem}[Corolary \ref{smoothness and irred property}]
Adopting the notations as above, then we have:\\
1. When $\sigma$ is not on any wall, $\mathfrak M_{\sigma}(n)$ is
either a smooth, irreducible variety of dimension $2n$ or empty.\\
2. Given $\sigma$ and $\sigma'$ not on any wall, $\mathfrak
M_{\sigma}(n)$ and $\mathfrak M_{\sigma'}(n)$ are birational to each
other when both are non-empty.\end{theorem}

To prove this theorem we study the GIT construction of $\mathfrak
M_{\sigma}(n)$ in detail, and control the dimension of the
exceptional locus for each birational map associated to wall
crossing. Then for each moduli space $\mathfrak M_{\sigma}(n)$, we
assign an ample line bundle on it. Applying the variation of
geometric invariant theory by \cite{DoHu} and \cite{Th}, we show
that a Bridgeland stability wall-crossing of $\mathfrak
M_{\sigma}(n)$ is the flip with respect to the line bundle. As a
result, the nested semicircular walls are one to one correspondence
to the stable base locus decomposition walls of the effective
divisor cone of Hilb$^n(\textbf P^2)$. In addition, given the the
location of the destabilizing wall, its corresponding base locus
decomposition wall is written out in an explicit way. Notice that in
certain cases, this correspondence has been established in a recent
paper \cite{CH} by Coskun and Huizenga via a different approach.

\begin{theorem}(Theorem \ref{main theorem in the
body}, Proposition \ref{proposition: the explicit formula of 1-1
walls}) For each semicircular actual wall on the second quadrant of
Bridgeland stability conditions plane we may assign a divisor
$\mathcal L_{\rho_{s,t,k}}$ up to a scalar. $\mathcal
L_{\rho_{s,t,k}}$ is on the stable base locus wall of Hilb$^nS$ and
this gives a one to one correspondence between the walls in the
stability plane and the stable base locus walls in the effective
divisor cone. In particular, the destabilizing semicircular wall on
the Bridgeland stability condition space with center
$-m-\frac{3}{2}$ corresponds to the base locus wall spanned by
divisor $mH-\frac{\Delta}{2}$. \label{main theorem in introducton}
\end{theorem}

Another important attempt in this paper is to extend this story to
the deformations of  Hilb$^n\textbf P^2$ by methods from
non-commutative algebraic geometry. Here we use the notion of
Sklyanin algebras $S$ $=$ Skl$(E,\sigma, \mathcal L)$, which are
non-commutative deformations of the homogeneous coordinate ring of
$\textbf P^2$. Such a Sklyanin algebra depends on a cubic curve $E$
on \textbf P$^2$, an automorphism $\sigma$ of $E$ and a degree $3$
line bundle. The foundation of such a non-commutative theory has
been set up in \cite{ATV1}, \cite{ATV2}, \cite{AV}, \cite{NS},
\cite{SV}. For these non-commutative $\textbf P^2$, we still have
$\mathfrak M^{ss}_{MG}(1,0,1-n)$, which turn out to be smooth
varieties (in the ordinary commutative sense!), and are in fact
deformations of Hilb$^n\textbf P^2$ by \cite{NS} and \cite{Hi}. We
will call these deformations of Hilb$^nS$.\\

In this paper we study the Bridgeland stability conditions of
$D^b(Coh(S))$, which are similar to that of $\textbf P^2$. In
particular, we have a similar chamber structure on the upper half
plane slice of the Bridgeland space, and the theorem above also
holds for $\mathfrak M_{\sigma}(n)$  associated to non-commutative
$\textbf P^2$. However, the behavior of wall-crossing over the
vertical wall is different in this case, and this changes the
correspondence between the chamber walls of stability space and the
base locus decomposition walls of the effective cone. In this case,
we have the following theorem:
\begin{theorem}(Theorem \ref{main theorem in the
body}, \ref{left half upper plane's main theorem in the body}) When
$n\geq 3$, for each semicircular actual wall on the Bridgeland
stability conditions plane we may assign a divisor $\mathcal
L_{\rho_{s,t,k}}$ up to a scalar. $\mathcal L_{\rho_{s,t,k}}$ is on
the stable base locus wall of Hilb$^nS$ and this gives a one to one
correspondence between the walls in the stability plane and the
stable base locus walls in the effective divisor cone. \label{main
theorem in introducton}
\end{theorem}
In addition this map is `monotone' in the sense that the two most
inner walls on the two quadrants correspond to the two edges of the
effective cone respectively. When one moves from inner semicircles
to the outside, the corresponding stable base locus wall moves in
one direction. This reveals a symmetric structure of the Mori
decomposition of the effective divisor cone of Hilb$^nS$. Notice
that, given any $n$, the destabilizing walls of Hilb$^n\textbf P^2$
and Hilb$^n$\textbf S are computable. Using the location of these
destabilizing walls, we can also compute the slopes of base locus
decomposition walls in the effective divisor cone. The cartoon of
divisor cones of Hilb$^n$\textbf P$^2$ and a generic Hilb$^n$\textbf
S are
shown below.\\\\\\

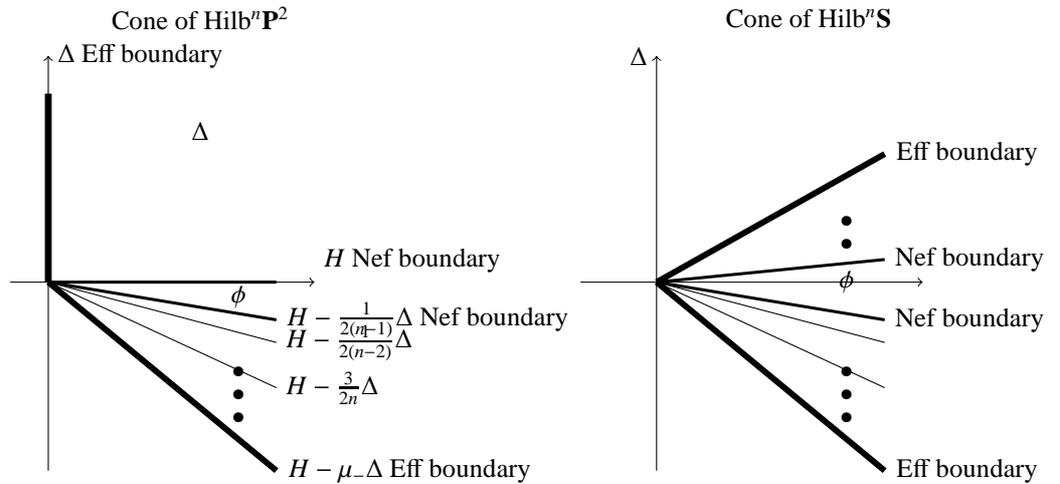
\begin{tikzpicture}[domain=0:1]
\draw (2,3.5) node {Cone of Hilb$^n$\textbf P$^2$};

\draw (10,3.5) node {Cone of Hilb$^n$\textbf S};

\draw[->] (-0.5,0) -- (3.5,0) node[above right] {$H$ Nef boundary};
\draw[->] (0,-2.5) -- (0,3) node[right] {$\Delta$ Eff boundary};
\draw [line width=0.4mm](0,0) --(3,-0.5) node[right]
{$H-\frac{1}{2(n-1)}\Delta$ Nef boundary}; \draw [line width=0.4mm]
(0,0) -- (3,0); \draw (0,0) --(3,-0.8) node[right]
{$H-\frac{1}{2(n-2)}\Delta$};  \draw (0,0) --(3,-1.4) node[right]
{$H-\frac{3}{2n}\Delta$};

\draw [line width=0.8mm] (0,0) --(3,-2.5)node[right]
{$H-\mu_-\Delta$ Eff boundary}; \draw [line width=0.9mm] (0,0)--
(0,2.5); \draw (2,2) node {$\Delta$}; \draw (2.5,-0.2) node
{$\phi$};\draw (2.5,-1.5) node {$\bullet$};\draw (2.5,-1.8) node
{$\bullet$};\draw (2.5,-1.2) node {$\bullet$};

\draw[->] (7,0) -- (11.5,0) node[above right] {}; \draw[->] (8,-2.5)
-- (8,3) node[left] {$\Delta$};

\draw [line width=0.4mm](8,0) --(11,-0.5) node[right] { Nef
boundary};

\draw (8,0) --(11,-0.8);
\draw (8,0) --(11,-1.4);

\draw [line width=0.8mm] (8,0) --(11,-2.5)node[right] { Eff
boundary};

\draw [line width=0.4mm] (8,0) --(11,0.3) node[right] {Nef
boundary};

\draw [line width=0.8mm] (8,0) --(11,1.7)node[right] { Eff
boundary}; \draw (10.5,.5) node {$\bullet$};\draw (10.5,.8) node
{$\bullet$};

\draw (10.5,0) node {$\phi$};\draw (10.5,-1.5) node
{$\bullet$};\draw (10.5,-1.8) node {$\bullet$};\draw (10.5,-1.2)
node {$\bullet$};
\end{tikzpicture}
\\

In the picture on the left, $\Delta$ is the exceptional divisor of
the Hilbert-Chow map to Sym$^n$\textbf P$^2$, and $H$ is the
pull-back of $\mathcal O(1)$ on Sym$^n$\textbf P$^2$. The picture on
the right is for Hilb$^n$\textbf S. Here $\Delta$ and $H$ are the
corresponding divisor classes under deformation. It is immediate
from the picture that Hilb$^n$\textbf S are
Fano and Hilb$^n$\textbf P$^2$ is log Fano.\\

After we obtained the results in this paper but before we finished
writing it, the paper \cite{CH} of Coskun and Huizenga appeared. In
\cite{CH}, the authors obtained the `correspondence of walls' result
for Hilb$^n\textbf P^2$ in certain cases. The paper \cite{CH} does
not treat the case of Hilb$^nS$, which is new in this paper. Also,
in \cite{CH} the author study the zero dimensional monomial
subschemes $Z$ of \textbf P$^2$, and when $\mathcal I_Z$ is
destabilized to get their result. our approach is quite different,
and the approach in \cite{CH} does not apply to the non-commutative
case, for example, only $n$-dimensional points in Hilb$^nS$
correspond to ideal sheaves. We show the smoothness and
irreducibility of each moduli space by showing some Ext$^2$
vanishing. These good properties allow one to apply the VGIT to get
the correspondence.\\\\
\textbf{Acknowledgments.} We are first of all greatly indebted to
Emanuele Macri, who offered tremendous assistance during the
preparation of this work. We are grateful to Izzet Coskun and Jack
Huizenga for helpful conversations. We also had useful discussions
with our advisors Herbert Clemens, Robert Lazarsfeld, Thomas Nevins
and Karen Smith, and we would like to thank all of them.

\section{Background Material}
\subsection{Review: Sklyanin Algebra and deformed Hilbert scheme of \textbf P$^2$}
We first recollect some definitions and properties about the
Sklyanin algebra from the noncommutative algebraic geometry,
further details are referred to \cite{NS} by Nevins and Stafford.\\

Given a smooth elliptic curve $\iota:$ $E\hookrightarrow$ \textbf
P$^2$ with corresponding line bundle $\mathcal L$ $=$
$\iota^*(\mathcal O_{\textbf P^2}(1))$ of degree $3$ and an
automorphism $\sigma$ $\in$ Aut$(E)$ which is a translation under
the group law. Denote the graph of $\sigma $ by $\Gamma_{\sigma}$
$\subset$ $E\times E$. Let $V$ $:=$ H$^0(E,\mathcal L)$, then we
have a $3$-dimensional space $\mathcal R(E,\sigma,\mathcal L)$:
\[
\mathcal R = \text{H}^0(E\times E,(\mathcal L \boxtimes\mathcal
L)(-\Gamma_{\sigma})\subset\text{ H}^0(E\times E,\mathcal L\boxtimes
\mathcal L)= V \otimes V.
\]
\begin{definition} The
\textit{$3$-dimensional Sklyanin algebra} is the algebra
\[
S=\text{Skl}(E,\mathcal L,\sigma)=T(V)/(\mathcal R),
\]
where $T(V)$ denotes the tensor algebra of $V$. \end{definition}
When $\sigma$ is the identical morphism, Skl$(E,\mathcal L,$ Id$)$
is just the commutative polynomial ring $\mathbb C[x,y,z]$. In
general, one may write Skl$(E,\mathcal L,\sigma)$ as a $\mathbb
C$-algebra with generators $x_1$, $x_2$, $x_3$ satisfying
relations:\[ax_ix_{i+1} +bx_{i+1}x_i+cx^2_{i+2}=0\text{    },
i=1,2,3\text{   mod } 3,\tag{$\triangle$}\label{relation of x,y,z}\]
where $a,b,c$ $\in$ $\mathbb C^*$ are scalars such that $(3abc)^3$
$\neq$ $(a^3+b^3+c^3)^3$.\\

$S=$ Skl$(E,\mathcal L, \sigma)$ is a connected graded algebra with
grading induced from $T(V)$. Write Mod-$S$ for the category of right
$S$-modules and Gr-$S$ for the category of graded right $S$-modules,
with homomorphisms Hom$_S(M,N)$ being graded homomorphisms of degree
zero. Given a graded module $M$ $=$ $\oplus_{i\in \mathbb Z} M_i$,
the shift $M(n)$ is the graded module with $M(n)_i$ $=$ $M_{i+n}$
for all $i$.
$S$ is strongly noetherian (Lemma 5.1 in \cite{NS}), we may write
gr-$S$ for the subcategory of noetherian objects in Gr-$S$. A module
$M$ $\in$ gr-$S$ is called right bounded if $M_i =0$ for $i\gg 0$.
The full Serre subcategory of gr-$S$ generated by the right bounded
modules is denoted by rb-$S$ with the quotient category qgr-$S$ $=$
gr-$S/$ rb-$S$. One has an adjoint pair $\pi:$ gr$-S
\leftrightarrows$ qgr-$S$ $:\Gamma^*$. Here $\pi$ is the natural
projection and $\Gamma^*$ is the `global section' functor. When $S$
is $\mathbb C[x,y,z]$, qgr-$S$ is isomorphic to the category of
coherent sheaves on \textbf{Proj}$\mathbb C[x,y,z]$. Due to this
reason, we call an
object $\mathcal M$ $\in$ qgr-$S$ as a sheaf on $S$.\\

A sheaf $\mathcal M$ on $S$ is called \textit{torsion} if each
element in $\Gamma^*(\mathcal M)$  is annihilated by a nonzero
element of $S$, respectively \textit{torsion-free} if no element is
so.  A torsion-free $\mathcal M$ has \textit{rank} $r$ if $\mathcal
M$ contains a direct sum of $r$, but not $r+1$, nonzero submodules.
The rank of a general sheaf $\mathcal M$ is defined to be the rank
of its torsion-free quotient part. We write
rk$(\mathcal M)$ for the rank of $\mathcal M$.\\

The \textit{first Chern class} $c_1(\mathcal M)$ is defined in
\cite{NS} Lemma 3.7 as the unique function $c_1:$ qgr-$S$
$\rightarrow$ $\mathbb Z$ with the following properties: additive on
short exact sequences; $c_1(\mathcal O_S(m))=m$ for all $m\in
\mathbb Z$. The Euler character on qgr-$S$ is defined as usual:
$\chi(\mathcal E,\mathcal F)$ $:=$ $\sum_i (-1)^i$ dim
Ext$^i(\mathcal E,\mathcal F)$ for $\mathcal E,$ $\mathcal F$ $\in$
qgr-$S$. $\chi(\mathcal F)$ $:=$ $\chi(\mathcal O,\mathcal F)$. The
\textit{Hilbert polynomial} of $\mathcal M$ is $p_{\mathcal M}(t): =
\chi(\mathcal M(t))$. The \textit{Mumford-Giesker slope }of
$\mathcal M$ is defined as $\mu^{MG}(\mathcal M)$ $:=$ $c_1(\mathcal
M)/$rk$(\mathcal M)$. A torsion-free sheaf $\mathcal M$ is called
\textit{Mumford-Giesker stable}, if for every non-zero proper
submodule $\mathcal F\subset \mathcal M$, one has rk($\mathcal
M)p_{\mathcal F}$ - rk($\mathcal F)p_{\mathcal M}$ $<$ $0$. Given a
torsion-free sheaf $\mathcal M$, it has a Harder-Narasimhan
filtration $0=\mathcal M_0\subset \mathcal M_1\dots \subset \mathcal
M_n=\mathcal M$ such that each quotient $\mathcal F_i =\mathcal
M_i/\mathcal M_{i-1}$ is Mumford-Giesker semistable with slope
$\mu^{GM}(\mathcal F_i)$ $>$ $\mu^{GM}(\mathcal F_{i+1})$. We write
$\mu^{GM}_+(\mathcal M)$ for $\mu^{GM}(\mathcal F_1)$, and
$\mu^{GM}_-(\mathcal M)$ for $\mu^{GM}(\mathcal F_n)$.
\begin{lemma}
Let $\mathcal M$ $\in$ qgr-$S$.\\
1. $c_1(\mathcal M(s))$ = $c_1(\mathcal M)$ $+$ $s\cdot$rk$(\mathcal
M)$ for any $s\in\mathbb Z$;\\
2. If $\mathcal M$ is torsion and non-zero, then $c_1(\mathcal M)$
$\geq$ $0$; if in addition $c_1(\mathcal M)$ $=$ $0$, then
$\chi(\mathcal M)$ $>$ $0$. \label{property of first chern class and
torsion sheaf}
\end{lemma}
\begin{proof}
Property 1 is the same as the second property of Lemma 3.7 in
\cite{NS}.\\

For property 2, let $\mathcal O(j)\rightarrow \mathcal M$ be a
non-zero morphism. By noetherian hypothesis on $\mathcal M$, the
descending chain of the quotient sheaves of $\mathcal M$ is finite.
By the additivity of $c_1$ and $\chi$, we may assume that $\mathcal
O(j)\rightarrow \mathcal M$ is surjective. To check that $c_1$ is
non-negative, by the first property, we may assume $j=0$. Let
$\mathcal I$ be the kernel of $\mathcal O\rightarrow \mathcal M$.
Denote $\Gamma^*(\mathcal I)$ by $I$. Write $c$ for $c_1(\mathcal
I)$. $\mathcal I(-c)$ is a rank $1$, normalized (i.e. $c_1(\mathcal
I(-c))=0$), torsion-free sheaf. By Proposition 5.6, Theorem 5.8 and
Lemma 6.4 in \cite{NS}, $\mathcal I(-c)$ is the homological sheaf
H$^0(\textbf K)$ of
\[ \textbf K: \mathcal O(-1)^{\oplus a}\rightarrow \mathcal O^{\oplus
2a+1}\rightarrow \mathcal O(1)^{\oplus a}\] at the middle term,
where $a$ is $1-\chi(\mathcal I(-c))$. Now for $n\gg 0$, recall $I$
$=$ $\oplus_{n\in\mathbb Z} I_n$, we have:
\begin{align*}
\text{dim}_{\mathbb C} I_n &  = (2a+1)\text{dim}_{\mathbb C} S(c)_n
- a\text{ dim}_{\mathbb C} S(c-1)_n - a\text{ dim}_{\mathbb C}S(c+1)_n\\
& =  (2a+1)\begin{pmatrix} n+c+2\\
2\end{pmatrix} - a(\begin{pmatrix}
n+c+1\\2\end{pmatrix}+\begin{pmatrix}
n+c+3\\2\end{pmatrix})\\
& =  \begin{pmatrix} n+c+2\\2\end{pmatrix} - a.
\end{align*}
Since $\mathcal I$ is a subsheaf of $\mathcal O$, dim$_{\mathbb C}
I_n$ $<$ dim$_{\mathbb C} S_n$ $=$ $\begin{pmatrix}
n+2\\2\end{pmatrix}$ for $n \gg 0$. We get $c\leq 0$, hence
$c_1(\mathcal M)\geq 0$.\\
When $rk(\mathcal M)$ $=$ $c_1(\mathcal M)$ $=$ $0$, by the formula
2 in Lemma 6.1 in \cite{NS}, the Hilbert polynomial $\chi(\mathcal
M(t))$ $=:$ $p_{\mathcal M}(t)$ $=$ $\chi(\mathcal M)$ is a
constant, we may also assume that $j=0$. Then $\mathcal I$ is
semistable and normalized, by Lemma 6.4 in \cite{NS}, $\chi(\mathcal
I)$ $\leq$ $1$ and the equality only holds when $\mathcal I$ $=$
$\mathcal O$.
\end{proof}
Let $D^b($qgr-$S)$ be the bounded derived category of qgr-$S$. We
rephrase one of the main results in \cite{NS}.
\begin{proposition}[Proposition 6.20 in \cite{NS}]
$D^b($qgr-$S)$ is generated by (i.e. the closure under that
extension and the
 homological shift of) $\mathcal O(k-1)$, $\mathcal O(k)$, $\mathcal
O(k+1)$ for any $k\in \mathbb Z$. \label{theorem:Db qgr-S is
generated by Ok's}
\end{proposition}
\begin{proof}
For any integer k, by the induction on $k$ and the exact sequence
\[
0\rightarrow \mathcal O(k)
\xrightarrow{\begin{bmatrix}z\\x\\y\end{bmatrix}} \mathcal
O(k+1)^{\oplus 3} \xrightarrow{\begin{bmatrix}ay & cx & bz\\bx & az
& cy\\cz & by & zx\end{bmatrix}}\mathcal O(k+2)^{\oplus 3}
\xrightarrow{\begin{bmatrix}x\\y\\z\end{bmatrix}^T}\mathcal
O(k+3)\rightarrow 0,
\]
where $a,b,c$ are coefficients in (\ref{relation of x,y,z}),
$\mathcal O(k)$ is in the closure. By Proposition 6.20 in \cite{NS},
all the Mumford-Giesker semi-stable sheaves are in the closure.
Since each torsion-free sheaf admits a finite Harder-Narasimhan
filtration, and each torsion sheaf is the cokernel of a morhism
between two torsion free sheaves, all sheaves are in the closure.
\end{proof}
As a consequence, invariants $\{$rank, first Chern class, Euler
character$\}$ generate the numerical Grothendieck group of
D$^b$(qgr-$S$). The importance of Sklyanin algebras is shown in the
following theorem in \cite{NS}. There the authors prove that
deformations of Hilb$^n$ \textbf P$^2$ can be constructed as the
moduli spaces of (semi)stable objects in qgr-$S$ with numerical
invariants $(1,0,1-n)$. As pointed out in \cite{Hi}, generically
each deformation of Hilb$^n\textbf P^2$  is constructed in this way.
\begin{theorem}[Theorem 8.11, 8.12  in \cite{NS}]
Let $\mathcal B$ be a smooth curve defined over $\mathbb C$ and let
$S_{\mathcal B}$ ($=$ $S_{\mathcal B}(E,\mathcal L,\sigma)$) $\in$
\underline{AS}$'_3$ be a flat family of algebras such that $S_p$ $=$
$\mathbb C[x,y,z]$ for some point $p$ $\in$ $\mathcal B$. Set $S$
$=$ $S_b$ for any point $b$ $\in$ $\mathcal B$. Then $\mathfrak
M^{ss}_{\mathcal B}(1,0,1-n)$ is smooth over $\mathcal B$, and
$\mathfrak M^{ss}_{\mathcal B}(1,0,1-n)$ $\otimes_{\mathcal B}$
$\mathbb C(b)$ $=$ $\mathfrak M^{ss}_{S_b}(1,0,1-n)$. Each
$\mathfrak M^{ss}_S(1,0,1-n)$ is a smooth, projective, fine moduli
space for equivalence classes of rank one torsion-free modules
$\mathcal M$ $\in$ qgr-$S$ with $c_1(\mathcal M)$ $=$ $0$ and
$\chi(\mathcal M)$ $=$ $1-n$. Each $\mathfrak M^{ss}_{S_b}(1,0,1-n)$
is a deformation of Hilb$^n$ \textbf P$^2$. \hfill $\square$
\end{theorem}
We will write $S^{[n]}$ or Hilb$^n S$ instead of $\mathfrak
M^{ss}_S(1,0,1-n)$ for short.
\begin{proposition}
The Picard number of Hilb$^nS$ is $2$.
\end{proposition}
Proof: By the formula on the second page of \cite{Nak} by Nakajima,
$b_2$ (Hilb$^nS$) $=$ $b_2$ (Hilb$^n$ \textbf P$^2$) $=$ $2$. Since
Hilb$^nS$ is projective, the Hodge numbers $h^{1,1} \geq 1$ and
$h^{0,2}$ $=$ $h^{2,0}$, one must have $h^{1,1}$ $=$ $2$.\hfill
$\square$
\subsection{Review: Bridgeland Stability Condition on D$^b($qgr-$S)$}
In this section, we briefly review the stability conditions on
derived categories. These notations are introduced in \cite{Br} by
Bridgeland. Let $N(S)$ be the numerical Grothendieck group of
$D^b($qgr-$S)$, i.e., the free abelian group generated by $r$,
$c_1$, and $\chi$.
\begin{definition}
A numerical stability condition on $D^b($qgr-$S)$ is a
pair\[(Z,\mathcal A), \mathcal A\subset D^b(\text{qgr-}S),\] where
$Z:$ $N(S)\rightarrow \mathbb C$ is a group homomorphism and
$\mathcal A$ is the heart of a bounded t-structure, such that the
following conditions hold.\\
1. For any non-zero $E$ $\in$ $\mathcal A$, we have \[Z(E)\in\{
re^{i\phi \pi}: r>0, \phi\in (0,1]\}.\] 2. Harder-Narasimhan
property: for any $E\in \mathcal A$, there is a filtration of finite
length in $\mathcal A$ \[0=E_0\subset E_1\subset \dots \subset
E_n=E\] such that each subquotient $F_i=E_i/E_{i-1}$ is
$Z$-semistable with arg $Z(F_i)>$ arg $Z(F_{i+1})$.
\label{definition: stability condition}
\end{definition}
Here an object $E$ $\in$ $\mathcal A$ is said to be $Z$-(semi)stable if for
any subobject $0$ $\neq$ $F$ $\subsetneq$ $E$ in $\mathcal A$ we
have
\begin{center} arg $Z(F)<(\leq)$ arg $Z(E)$.\end{center}
The group homomorphism $Z$ is called the central charge of the
stability condition.
The rest of this section is devoted to the construction of numerical
stability conditions on qgr-$S$. First we recall the notion of
torsion pairs, which is essential to constructing t-structures. A
pair of full subcategories $(\mathcal F,\mathcal T)$ of qgr-$S$ is
called a torsion pair if it satisfies the following two conditions.
\begin{enumerate}
\item For all $F$ $\in$ ob$\mathcal F$ and $T$ $\in$ ob$\mathcal
T$, we have Hom$(T,F)$ $=$ $0$.

\item Each sheaf $E$ in qgr-$S$ fits in a short exact
sequence:\begin{center} $0$ $ \rightarrow T$ $\rightarrow E$
$\rightarrow F$ $\rightarrow 0$,
\end{center}
where $T\in$  ob$\mathcal T$, $F\in$ ob$\mathcal F$. In addition,
the extension class is uniquely determined up to isomorphism.
\end{enumerate}

A torsion pair defines a t-structure on $D^b($qgr-$S)$ by:
\begin{center} $\mathcal D^{\geq 0}$ $=$ $\{$ $C^{\bullet}|$
H$^{-1}(C^{\bullet})$ $\in$ $\mathcal F$ and H$^i(C^{\bullet})=0$
for
$i<-1$ $\}$,\\
$\mathcal D^{\leq 0}$ $=$ $\{$ $C^{\bullet}|$ H$^{0}(C^{\bullet})$
$\in$ $\mathcal T$ and H$^i(C^{\bullet})=0$ for $i>0$
$\}$.\end{center} As in the \textbf P$^2$ case, given $s$ $\in$
$\mathbb R$, we can define the full subcategories $\mathcal T_s$ and
$\mathcal F_s$ of qgr-$S$ as:
\begin{flushleft}
 $T\in\mathcal
 T_s$ if $T$ is torsion or $\mu^{MG}_{-}(T')$ $>$ $s$, where $T'$ is
  the torsion-free quotient of $T$;\\
$F\in\mathcal F_s$ if $F$ is torsion-free and $\mu^{MG}_+(F)$ $\leq$
$s$.
\end{flushleft}
By Lemma 6.1 in \cite{Br2}, $(\mathcal F_s,\mathcal T_s)$ is a
torsion pair. Let $\mathcal A_s$ be the heart of the t-structure
determined by the torsion pair $(\mathcal F_s,\mathcal T_s)$, we may
define a central charge $Z_{t,s}$ $=$ $-d_{t,s}$ $+$ $ir_{t,s}$
depending on a parameter $t>0$ by:
\begin{flushleft}
$r_{t,s}(E)$ $:=$ $(c_1-rs)t$;\\
$d_{t,s}(E)$ $:=$ $-rt^2/2$ $+$ $(s^2/2-1)r$ $-$ $(3/2+s)c_1$ $+$
$\chi$.
\end{flushleft}
We may write $\mu_{s,t}$ = $d_{s,t}/r_{s,t}$ as the `slope' of an
object in $\mathcal A_s$.
\begin{proposition}
$(Z_{s,t},\mathcal A_s)$ is a stability condition on D$^b($qgr-$S)$.
\end{proposition}
\begin{proof} Any object $E$ $\in$ $\mathcal A_s$ fits in an exact
sequence:
\[
0\rightarrow \text{H}^{-1}(E)[1]\rightarrow E\rightarrow
\text{H}^0(E)\rightarrow 0.
\]
in order to check the Property 1 of the central charge in Definition
\ref{definition: stability condition}, we only need to check
arg$(Z(E))$ $\in$ $(0,\pi]$ for the following cases: 1. $E$ is a
torsion sheaf. 2. $E$ is a Mumford-Giesker stable sheaf with
$\mu(E)$ $>$ $s$. 3. $E[-1]$ is a Mumford-Giesker stable sheaf with
$\mu(E[-1])$ $\leq$ $s$.\\

Case 1 is due to Lemma \ref{property of first chern class and
torsion sheaf}. Case 2 is clear since $r_t(E)$ is greater than $0$.
In case 3, we may assume $c_1(E[-1])=r(E[-1])s$, then $d_t(E)$ $=$
$rt^2/2$ $+$ $r$ $+$ $3c_1/2$ $-$ $\chi$ $+$ $c_1^2/2r$ $\geq$
$rt^2/2$ $+$ $r^2$ $-$ $1$ $>$ $0$, where $r,c_1,\chi$ stands for
$r(E[-1])$, $c(E[-1])$, $\chi(E[-1])$ respectively. The first
inequality is due to Corollary 6.2 and Proposition 2.4 in \cite{NS}:
$2\chi r$ $-$ $r^2$ $-$ $3rc_1$ $-$ $c^2_1$ $=$ $\chi(E[-1],E[-1])$
$\leq$ $1$ $+$ ext$^2(E[-1],E[-1])$ $=$ $1$ $+$
h$^0(E[-1],E[-1](-3))$ $=$
$1$.\\

When $t$ and $s$ are both rational numbers, the Harder-Narashimhan
property is similarly proved  as that in Proposition 7.1 in
\cite{Br2} by Bridgeland. In the general case, we need the following
lemma to check the descending chain stable property:

\begin{lemma}
Given a stability condition $( Z_{s,t}, \mathcal A_s)$ and two
positive numbers $M_1, M_2$. Then the set
\[\big( (-\infty,M_1]\times [0,M_2] \big)\bigcap \{(-d_{s,t}(\mathcal
F),r_{s,t}(\mathcal F))|\mathcal F \in \mathcal A_s \text{ and is a
torsion-free sheaf}\}\] is finite. \label{only finite many r,d in
shaded area}
\end{lemma}
\begin{proof}
First, we show this holds for all Mumford-Giesker stable sheaves.
Write $\chi(\mathcal F)$, $r(\mathcal F)$, $c_1(\mathcal F)$ as
$\chi$, $r$, $c_1$ for short. By the inequality $2\chi
r-r^2-3rc_1-c^2_1\leq 1$, we have
\[\chi\leq \frac{1}{2r} (1+r^2+3rc_1+c^2_1).\]
Plug this into the formula of $d_{s,t}$, we have
\[ -d_{s,t}\geq \frac{rt}{2}-\frac{(c-rs)^2}{2r}.\]
Since $c_1-rs\geq 0$ and $M_1\geq -d_{s,t}$,
 $r$ is bounded. As $c_1-rs\in[0,M_2]$, $c_1$
is bounded. Now by the inequality of $\chi$ and the formula of
$d_{s,t}$, we have \[
-M_1+\frac{rt^2}{2}-(\frac{s^2}{2}-1)r-(\frac{3}{2}+s)c_1\leq
\chi\leq \frac{1}{2r} (1+r^2+3rc_1+c^2_1).\] Hence $\chi$ is
bounded. The set
\[((-\infty,M_1]\times [0,M_2])\bigcap \{(-d_{s,t}(\mathcal
F),r_{s,t}(\mathcal F))|\mathcal F \in \mathcal A_s \text{ and is a
torsion-free MG stable sheaf}\}\] is finite.\\

Next, we show this holds for torsion free $F\in \mathcal A_s$. By
the finiteness result above, we may define
\[D:=\min\{-d_{s,t}(\mathcal F)|\mathcal F \in \mathcal A_s \text{
and is a torsion-free MG stable sheaf}\};\]
\[R:=\min\{r_{s,t}(\mathcal
F)|\mathcal F \in \mathcal A_s \text{ and is a torsion-free MG
stable sheaf}\}.\] Now given a torsion free sheaf $\mathcal G\in
\mathcal A_s$, if $\mathcal G$ has an MG-factor $\mathcal G_i$ such
that $-d_{s,t}(\mathcal G_i)> \frac{D^2}{R}+M_1$, then
\[
-d_{s,t}(\mathcal G)>(\frac{D^2}{R}+M_1)+ D\cdot(-\frac{D}{R})
=M_1.
\]
Therefore \begin{align*} &\{(-d_{s,t}(\mathcal F),r_{s,t}(\mathcal
F))|\mathcal F \in \mathcal A_s \text{ and is a torsion-free
sheaf}\}\bigcap\big((-\infty,M_1]\times [0,M_2]\big)\\ \subset &
\sum \{(-d_{s,t}(\mathcal G),r_{s,t}(\mathcal G))|\mathcal G \in
\mathcal A_s \text{ and is a torsion-free MG stable
sheaf}\}\bigcap\big((-\infty, \frac{D^2}{R} + M_1]\times
[0,M_2]\big),
\end{align*}
and is finite.
\end{proof}

Now we may check the descending chain stable property. Suppose the
property does not hold, we have an object $E$ in $\mathcal A_{s,t}$
that has an infinite descending chain:
\[\dots \subset E_{i+1} \subset E_{i} \dots \subset E_{1}\subset E_{0} =E\]
with increasing slopes $\mu_{s,t}(E_{i+1})$ $>$ $\mu_{s,t}(E_i)$ for
all $i$. There are short exact sequences in $\mathcal A_{s,t}$:
$0\rightarrow E_{i+1}\rightarrow E_i\rightarrow F_i\rightarrow 0$
for $i\geq 0$. By taking the cohomology of sheaves, we have:
H$^{-1}(E_{i+1})\subset$ H$^{-1}(E_i)$. We may assume the rank of
H$^{-1}(E_{i})$ is constant. Now the cokernel H$^{-1}(E_{i})/$
H$^{-1}(E_{i+1})$ is torsion, and H$^{-1}(F_i)$ is torsion-free, we
have H$^{-1}(E_{i})$ $\simeq$ H$^{-1}(E_{i+1})$.\\

Let $T_i$ and $G_i$ be the torsion subsheaf and torsion-free
quotient of H$^0(E_i)$ respectively. Since we have  the exact
sequence
\[
0\rightarrow H^{-1}(F_i)\rightarrow H^0(E_{i+1})\rightarrow
H^0(E_i)\rightarrow H^0(F_i)\rightarrow 0,
 \]
 and H$^{-1}(F_i)$ is torsion-free, $T_{i+1}$ is a
subsheaf of $T_i$. We may assume $c_1(T_i)$ is constant, then
$\chi(T_i)$ is non-increasing.\\

Now we have
\[
\big(-d_{s,t}(H^{-1}(E_i)[-1])-d_{s,t}(T_i),r_{s,t}(H^{-1}(E_i)[-1])+r_{s,t}(T_i)\big)
\geq
\big(-d_{s,t}(H^{-1}(E_0)[-1])-d_{s,t}(T_0),r_{s,t}(H^{-1}(E_0)[-1])+r_{s,t}(T_0)\big).
\]
Since the slope is increasing, we also have
\[\big(-d_{s,t}(E_i),r_{s,t}(E_i)\big)\leq
\big(\max\{-d_{s,t}(E_0),0\},r_{s,t}(E_0)\big).\] Subtracting the
first from the second one, we have
\[\big(-d_{s,t}(G_i),r_{s,t}(G_i)\big)\leq
\big(\max\{-d_{s,t}(G_0),d_{s,t}(H^{-1}(E_0)[-1]+d_{s,t}(T_0),0\},r_{s,t}(G_0)\big).\]Now
applying the lemma to $G_i$, combining with the results on
$H^{-1}(E_i)[-1]$ and $T_i$, the set of possible values of
$\big(-d_{s,t}(E_i),r(E_i)\big)$ is finite, so we may get the
stableness of the descending chain directly as that in $s,t$
rational case.\\

The ascending chain property can be similarly proved, where one
applies the lemma to $H^{-1}(E_i)[-1]$ and the area $(-M_1,
+\infty)\times(0,M_2)$ to get the finiteness. Then the argument is
the same as that in $s,t$ rational case, and the details are left to
the readers.
\end{proof}

Let $\mathcal A(k)$ be the extension closure of $\mathcal
O(k-1)[2]$, $\mathcal O(k)[1]$ and $\mathcal O(k+1)$. Since
$<\mathcal O(k-1),\mathcal O(k),\mathcal O(k+1)>$ is a full strong
exceptional collection by Proposition \ref{theorem:Db qgr-S is
generated by Ok's}, $\mathcal A(k)$ is the heart of a t-structure of
D$^b($ qgr-$S)$, see Lemma 3.16 \cite{Ma}. Objects in $\mathcal
A(k)$ are of the form:
\[\mathcal O(k-1)\otimes \mathbb C^{n_{-1}}\rightarrow \mathcal O(k)\otimes \mathbb C^{n_0}\rightarrow \mathcal O(k
+1)\otimes \mathbb C^{n_1},\] where $n_{-1}$, $n_0$, $n_1$ are some
non-negative integers. We write $\overrightarrow{n}$ $=$
$(n_{-1},n_0,n_1)$ and call it the type of the object. One may
construct a central charge $Z$ for $\mathcal A(k)$ by letting
$Z(\mathcal O(k-1)[2])=z_{-1}$, $Z(\mathcal O(k)[1])=z_{0}$ and
$Z(\mathcal O(k+1))=z_{1}$ for any collection of complex numbers
$z_i$'s on the upper half plane: $\{ re^{i\phi \pi}: r>0, \phi\in
(0,1]\}$. $(Z,\mathcal A(k))$ is a stability condition on D$^b($
qgr-$S$).
\section{Destabilizing Wall}
The destabilizing walls on the $(s,t)$-plane of stability conditions
of D$^b($ coh\textbf P$^2)$ are discussed in \cite{ABCH} Section 6.
In the D$^b($ qgr-$S)$ case, the behavior of the walls
is similar to that of the \textbf P$^2$ case.\\

The potential wall associated to a pair of invariants $(r,c,\chi)$
and $(r',c',\chi')$ on qgr-$S$ is the following subset of the
upper-half $(s,t)$-plane:
\[
W_{(r,c,\chi),(r',c',\chi')}:=\{(s,t)|\mu_{s,t}(r,c,\chi)=\mu_{s,t}(r',c',\chi')\}.
\]
More explicitly, the wall is given by:
\[
W_{(r,c,\chi),(r',c',\chi')}=\{(s,t)| \frac{1}{2}(c_1r'-c_1'r)(t^2+s^2)+(\chi'r-\chi r'+\frac{3}{2}r'c_1-\frac{3}{2}rc_1')s + (c_1r'-c_1'r+\chi c_1'-\chi'c_1) =0\}.
\]

Let $W^{\text{potential}}_{r,c,\chi}$ $:=$ $\bigcup_{(r',c',\chi')}
W_{(r,c,\chi),(r',c',\chi')}$. In the Hilbert scheme case, where
$(r,c_1,\chi)$ $=$ $(1,0,1-n)$ (respectively $(-1,0,n-1)$ when
$s\geq 0$), the potential walls form the set
$\{(s,t)|-\frac{c_1'}{2}(s^2+t^2) +(\chi'-(1-n)r'-\frac{3}{2}c_1')s
-n c_1'=0\}$. When $c_1'$ $=$ $0$, the wall is the $t$-axis. When
$c_1'$ $\neq$ $0$, these are nested semicircles with center $x$ $=$
$(\chi'-(1-n)r'-\frac{3}{2}c_1')/c_1'$ and radius $Rad$ $=$
$\sqrt{x^2-2n}$. It is not hard to see that
$W^{\text{potential}}_{1,0,1-n}$ on the second quadrant and
$W^{\text{potential}}_{-1,0,n-1}$ on the first quadrant are nested
semicircles.\\

We define $W^{\text{actual}}_{r,c,\chi}$ as
\begin{center}
$\{(s,t)|\exists$ $\mathcal F$ with invariant $(r,c_1,\chi)$, which
is strictly semistable and locally stable under $(Z_{s,t},\mathcal
A_s)\}$.
\end{center}
Here by `$\mathcal F$ is
locally stable under $(Z_{s,t},\mathcal A_s)$', we mean that for any
$\delta
>0$, there is $(s',t')$ $\in$ $B_\delta(s,t)$ such that $\mathcal
F$ is stable under $(Z_{s',t'},\mathcal A_s')$. By definition,
$W^{\text{actual}}_{r,c,\chi}$ $\subset$
$W^{\text{potential}}_{r,c,\chi}$. On the second quadrant,
$W^{\text{actual}}_{1,0,1-n}$ is formed by nested semicircular
walls. We call such a wall in $W^{\text{actual}}_{1,0,1-n}$ as an
actual wall.
\begin{lemma}
For any $k\in \mathbb Z$, $\mathcal O(k)$ (resp. $\mathcal O(k)[1]$)
is a stable object under stability condition $(Z_{s,t}, \mathcal
A_s)$ for $s<k$ (resp. $s\geq k$).\label{O is stable}
\end{lemma}
\begin{proof}
Since $Z_{s,t}(E)$ $=$ $Z_{s+k,t}(E(k))$, we may assume $k$ is $0$.
Suppose $\mathcal O$ is not stable for some $(Z_{s,t},\mathcal
A_s)$, with $s<0$, then there exists $E$ destabilizing $\mathcal O$
under $(Z_{s,t},\mathcal A_s)$. We may assume $(s,t)$ is on a
potential wall of $W^{\text{potential}}_{1,0,1}$ which is a
semicircle with right corner at the origin. The exact sequence:
\[
0\rightarrow H^{-1}(E)\rightarrow H^{-1}(\mathcal O)\rightarrow
H^{-1}(\mathcal O/E) \rightarrow H^{0}(E)\rightarrow H^{0}(\mathcal
O)\rightarrow H^{0}(\mathcal O/E)\rightarrow 0
\]
implies that H$^{-1}(E)$ is $0$, we may write $E$ for H$^0(E)$ for
short. In addition, H$^{0}(\mathcal O/E)$ has rank $0$ (otherwise
the morphism $H^{0}(E) \rightarrow$ H$^{0}(\mathcal O)$ is $0$),
hence it has non-negative
$c_1$. This implies $\mu^{GM}_-(E)$ $<$ $0$.\\

Let $\mu^{GM}_-(E)$ $=$ $s_0$ $<$ $0$, we may move $(s,t)$ along the
semicircle to the right, when $s<s_0$, $E$ still destabilizes
$\mathcal O$ since $E$ and $\mathcal O$ are in $\mathcal A_s$. Write
$E$ as $0\rightarrow E_+\rightarrow E\rightarrow E_-\rightarrow 0$,
where $E_-$ stands for the MG semistable factor with slope $s_0$.
When $s$ tends to $s_0$, $r_{s,t}(E_{-})$ will tend to $0$. As
$d_{s,t}(E_{-})>l>0$ for some constant $l$, $E_+$ destabilizes
$\mathcal O$ under $(Z_{s_0,t},\mathcal A_{s_0})$. By repeating this
procedure, we get an $E'$ which destabilize $\mathcal O$ and has
$\mu_-(E')\geq 0$, this is a contradiction.
\\

For $\mathcal O[1]$ and $s> 0$ case, we get $\mu^{GM}_+($
H$^{-1}(\mathcal O/E))$ $>$ $0$, then the same argument also works.
When $s=0$, the only exceptional case is that both H$^{-1}(E)$ and
H$^{-1}(\mathcal O/E)$ has $c_1$ $=$ $0$ and H$^0(E)$ is torsion.
But that cannot happen since H$^{-1}(E)$ and H$^{-1}(\mathcal O/E)$
are torsion free and the rank of $\mathcal O$ is $1$.
\end{proof}

The $\widetilde{GL(2,\mathbb R)^+ }$ acts on the space of stability
condition by the SL$(2,\mathbb R)$-action on the central charge and
the shift on the heart structure. In particular, an element $\phi$
in the subgroup $\mathbb R$ acts on $(Z,\mathcal A)$ as follow: if
$\phi$ is an integer, then $\phi$ $\circ$ $(Z,\mathcal A)$ $=$
$(Z[\phi],\mathcal A[\phi])$ with:
\begin{center}
$\mathcal A[\phi]$ $=$ $\{A[\phi]|$ $A$ is an object of $\mathcal
A$\} and $Z[i](A)$ $:=$ $(-1)^iZ(A)$.
\end{center}
If $0<\phi<1$, then $\mathcal A[\phi]$ $:=$ $\langle \mathcal
T_{\phi},\mathcal F_{\phi}[1]\rangle$ and $Z[\phi](E)$ $:=$
$e^{-i\pi \phi} Z(E)$ with:
\begin{center}
$\mathcal T_{\phi}$ $=$ $\langle T\in \mathcal A|$ $T$ is stable
with arg$(Z(T))$ $>$ $\phi\pi\rangle$;

$\mathcal F_{\phi}$ $=$ $\langle F\in \mathcal A|$ $F$ is stable
with arg$(Z(F))$ $\leq$ $\phi\pi\rangle$.
\end{center}
The moduli spaces of stable objects are unaffected under this
$\mathbb R$-action.
\begin{proposition}[Proposition 7.5 in \cite{ABCH}]Let $k$ be an
integer. If a pair of real numbers $(s,t)$ satisfies
\[(s-k)^2 +t^2 < 1,\]
then there is $\phi_{s,t,k}\in\mathbb R$ (not canonically defined),
such that under its action $(Z_{s,t}[\phi_{s,t,k}],\mathcal
A_s[\phi_{s,t,k}])$ is identified with $(Z,\mathcal A(k))$ for
suitable choice of central charge $(z_{-1},z_0,z_1)$ for
$Z$.\label{Ast[phi] = Ak}
\end{proposition}
\begin{proof}
By Lemma \ref{O is stable}, the rest is the same as Prop 7.5
in \cite{ABCH}.
\end{proof}

We call such a semidisc a \emph{quiver
region}.\\

Consider a central charge $(z_{-1},z_0,z_1)$ of $\mathcal
A(k)$:\begin{center} $(z_{-1},z_0,z_1)$ $=:$ $\overrightarrow{z}$
$=$ $\overrightarrow{a} + i\overrightarrow{b}$,
\end{center} where $\overrightarrow{a}$ and $\overrightarrow{b}$ are
real vectors. Fix three non-negative integers $(n_{-1},n_0,n_1)$ $=$
$\overrightarrow{n}$, and let
\[\overrightarrow{\rho}=-\overrightarrow{a}+\overrightarrow{b}(\frac{\overrightarrow{n}\cdot\overrightarrow{a}}{\overrightarrow{n}\cdot\overrightarrow{b}})
\]
, then $\overrightarrow{n}\cdot \overrightarrow{\rho}$ $=$ $0$. An
object $E$ in $\mathcal A(k)$ with type $\overrightarrow{n}$ is
stable (semistable) with respect to the central charge
$\overrightarrow{z}$ if and only if for any proper subobject $E'$
with type $\overrightarrow{n}'$ one has:\begin{center}
$\overrightarrow{n}'\cdot\overrightarrow{\rho}$ $<$ $0$ ($\leq
0$).\end{center}
\begin{remark}\textit{
$\overrightarrow{\rho}$ does not change under the rotation of
$\overrightarrow{z}$, hence it does not depend on the choice of
$\phi_{s,t,k}$ in \ref{Ast[phi] = Ak}. The explicit formula of
$\overrightarrow{\rho}$ is given at \ref{character rho of s,t}.
}
\end{remark}
In particular, we will use the following computation in the $k=0$
 case of in the first statement of Proposition \ref{finite many wall,
t>0 is HilbS}.
\begin{example} For $\mathcal A(0)$ and $s<0$, let $\overrightarrow
{n}$ be $(n,2n+1,n)$, the  $\overrightarrow{\rho}_{s,t}$ is given
as:
\[\frac{ts}{\frac{t^2}{2} +n-
\frac{s^2}{2}}(\frac{(1+s)^2}{2}-\frac{t^2}{2},
\frac{t^2}{2}-\frac{s^2}{2},\frac{(1-s)^2}{2}-\frac{t^2}{2}) +
t(1+s,-s,s-1).\] \label{character rho of s,t for k=0}
\end{example}

Consider the space of characters $\overrightarrow{\rho}$, since the
subobjects of $E$ have only finitely many possible numerical types,
there are finitely many walls (lines in this case) on which an
object $E$ with type $\overrightarrow{n}$ could be semistable but
nonstable with respect to $\overrightarrow{\rho}$. These walls
divide the space into chambers. In each chamber, the moduli space of
stable objects remains the same, so one may choose an integral
vector $\overrightarrow{\rho}$ as a representative in the chamber.
By Proposition 3.1 in \cite{Ki} by King, the moduli space of
(semi)stable object with respect to central charge $Z$ consists the
$\overrightarrow{\rho}$-(semi)stable points under the $G$-action. As
explained in the Chapter 2.2 in \cite{Gi} by Ginzburg, the moduli
space of $Z$-semistable objects is constructed as a GIT quotient:
\begin{center}
Proj$\bigoplus_{n\geq 0} \mathbb C[X]^{G,\overrightarrow{\rho}^n}$.
\end{center}
Here $X$ is the affine closed subscheme of Hom$(\mathbb
C^{n_{-1}}\otimes \mathcal O(k-1), \mathbb C^{n_0} \otimes \mathcal
O (k))$ $\times$ Hom$(\mathbb C^{n_{0}}\otimes \mathcal O(k),
\mathbb C^{n_1} \otimes \mathcal O (k+1))$ consisted of the complex.
$G$ is the reductive group GL$(n_{-1},\mathbb C)$ $\times$
GL$(n_{0},\mathbb C)$ $\times$ GL$(n_{1},\mathbb C)/\mathbb
C^\times$ and $\overrightarrow{\rho}$ is the character $($
det$^{\rho_{-1}}$, det$^{\rho_{0}}$, det$^{\rho_{1}}$) of $G$. This
character is well-defined since
$\overrightarrow{\rho}\cdot\overrightarrow {n}$ is $0$. When
$\overrightarrow {n}$ is primitive (i.e. gcd$(n_{-1},n_0,n_1)$ $=$
$1$), $G$ acts freely on the stable points on $X$. We write
$\mathfrak M^{\overrightarrow{\rho},ss}(\overrightarrow {n})$ $:=$
$X$/$_{\overrightarrow{\rho}}G$ as the moduli space of semistable
objects in $\mathcal A(k)$ with type $\overrightarrow {n}$ and
character $\overrightarrow{\rho}$.
\begin{proposition}
1. Given $n>0$, for any $s<0$ and $t\gg 1$, the moduli space of
stable objects with invariants $(r,c_1,\chi)$ $=$ $(1,0,1-n)$ under
$(Z_{s,t},\mathcal A_s)$ is the same as that in Mumford-Giesker
sense, i.e., the moduli
space is the deformed Hilbert scheme Hilb$^nS$.\\
2. There are only finitely many actual destabilizing walls for Hilb$^nS$.
\label{finite many wall, t>0 is HilbS}
\end{proposition}
\begin{proof} Let $\mathcal I$ be a torsion free sheaf with  $(r,c_1,\chi)$
$=$ $(1,0,1-n)$. When $k =0$, $\mathcal I[1]$  appears in $\mathcal
A(0)$ with type $\overrightarrow {n}$ $=$ $(n,2n+1,n)$. By
Proposition 7.7 and Proposition 6.20 in \cite{NS}, let
$\overrightarrow{\rho}$ be $((2n+1)(m-1),n,-(2n+1)m)$, then for all
$m\gg 1$, $X^{ss,\overrightarrow{\rho}}$ consists of complexes which
are quasi-isomorphic to $\mathcal I[1]$ for some torsion-free
$\mathcal I$ with invariants $(1,0,1-n)$. Now by the formula in
Example \ref{character rho of s,t for k=0}, there is an open area
$A$ in the region $\{(s,t)| s^2+t^2<1, s<0\}$ with boundary
containing $(0,0<t<1)$ such that the stable objects with invariants
$(1,0,1-n)$ under $(Z_{s,t}, \mathcal A_s)$
are the same as those in the Mumford-Giesker sense. \\

Note that when $s<0$, $W^{\text{actual}}_{1,0,1-n}$ consists of
semicircles with center at $x$ and radius $\sqrt{x^2-2n}$. Since
$x+\sqrt{x^2-2n}$ is decreasing when $x<\sqrt{2n}$, these
semicircles are nested with right boundary in the region
$(-\sqrt{2n}<s<0,t=0)$, hence all the actual destabilizing walls are
nested semicircles and each of them corresponds a wall in $\mathcal
A(k)$ for some $0<k<\sqrt{2n}$. This tells the finiteness of the
actual walls, and in the region outside the first wall, the stable
objects with invariants $(1,0,1-n)$ are
the same as those in the area $A$. \\

When $s>0$, the same argument works for the second statement.
\end{proof}
\section{Wall Crossing via GIT}
\subsection{Properties of stable objects in $\mathcal A(k)$}
The goal of this section  is to show Corollary \ref{smoothness and
irred property}: the moduli space $\mathfrak M^{\overrightarrow
{\rho},s}(\overrightarrow{n})$ is smooth and irreducible for generic
$\overrightarrow{\rho}$. The next two lemmas are about the vanishing
property of some Ext$^2$'s.
\begin{lemma}
Let $\mathcal F$ be a stable object in $(\mathcal A_s,Z_{s_0,t_0})$
(for some$ s_0<0$) with invariant $(r,c_1,\chi)$ $=$ $(1,0,1-n)$.
Then we have 
\begin{center}
\emph{Hom}$(\mathcal F,\mathcal F[2]) = 0. $\end{center}
\label{smoothness lemma}
\end{lemma}
\begin{proof}
Given a point $(\tilde{s},\tilde{t})$ on the second quadrant, we
denote $W_{(\tilde{s},\tilde{t})}$ as  the unique semicircle with
central at $x$ and radius $\sqrt{x^2-2n}$ that across
$(\tilde{s},\tilde{t})$, i.e. the potential semicircular wall for
the invariant $(1,0, 1-n)$ across
$(\tilde{s},\tilde{t})$.\\\\
\textbf{Case I}: $W_{(s_0,t_0)}$ has radius greater than
$\frac{3}{2}$. The actual destabilizing walls of $\mathcal F$ are
nested, so $\mathcal F$ is a stable object under $(\mathcal
A_s,Z_{s,t})$ for all $(s,t)$ $\in$ $W_{(s_0,t_0)}$. $\mathcal
F(-3)$ is a stable object under $(\mathcal A_{s-3},Z_{s-3,t})$, for
any $(s,t)$ $\in$ $W_{(s_0,t_0)}$. These points form the semicircle
$W_{(s_0,t_0)}$ $-$ $(3,0)$. Since the radius of the circle is
greater than $\frac{3}{2}$, these two semicircles intersect at
$(s_1,t_1)$. In $\mathcal A_{s_1}$, under the central charge
$Z_{s_1,t_1}$, the slope of $\mathcal F$ is
$-\frac{s_1}{2}+\frac{t_1^2+2n}{2s_1}$. Because $s_1$ is less than
$-3$,  $-\frac{s_1}{2}+\frac{t_1^2+2n}{2s_1}$ is greater than the
slope of $\mathcal F(-3)$, which is
$-\frac{s_1+3}{2}+\frac{t_1^2+2n}{2s_1+3}$. Thus Hom$(\mathcal
F,\mathcal F(-3))$ $=$ $0$, and Hom($\mathcal
F,\mathcal F[2]$) $=$ $0$ by Serre duality.\\\\
\textbf{Case II}: $W_{(s_0,t_0)}$ has radius equal to or less than
$\frac{3}{2}$. Let $k$ be the positive integer such that
\[ (k+1)(k+2)/2\leq n < (k+2)(k+3)/2.\]
The semicircle $W_{(-k-1,0)}$ has radius not less than
$\frac{1}{2}$, and by Lemma \ref{lemma: existence of the last wall}
after this wall, there is no stable object with invariant
$(1,0,1-n)$. The radius of $W_{(-k,0)}$ is greater than
$\frac{3}{2}$, hence the right edge of $W_{(s_0,t_0)}$ falls into
the interval $(-k-1,-k)$. Therefore $\mathcal F[1]$ is an object in
$\mathcal A(-k)$, and
$\mathcal F(-3)[1]$ is an object in $\mathcal A(-k-3)$.\\

On the other hand, $W_{(s_0,t_0)}$ is larger than $W_{(-k-1,0)}$,
hence its left edge is less than $-k-2$. Since its radius is not
greater than $\frac{3}{2}$, its left edge is greater than $-k-4$.
Combining these two observations, the left edge of $W_{(s_0,t_0)}$
falls into the $\mathcal A(-k-3)$ quiver region. Therefore $\mathcal
F$ is an object in $\mathcal A(-k-3)$. We have
\begin{center}
Hom ($\mathcal F,\mathcal F(-3))$ $=$ Hom $(\mathcal F,(\mathcal
F(-3)[1])[-1])$ $=$ $0$,
\end{center}
where the last equality is because of both $\mathcal F$ and
$\mathcal F(-3)[1]$ are in a same heart $\mathcal A(k+3)$. By Serre
duality, Hom($\mathcal F,\mathcal F[2])$ $=$ $0$.
\end{proof}

On each chamber wall all S-equivalent semistable objects (i.e.,
their stable factors are the same after rearrangement) are
contracted to one point. Let $\mathcal F$ be a locally stable object
with invariants $(1,0,1-n)$ at $(s_0,t_0)$. Suppose it is
destabilized at this point which lies on an semicircular actual wall
$W_{(s_0,t_0)}$. Then $\mathcal F$ has a filtration in $\mathcal
A_{s_0}$:
\[\mathcal F = \mathcal F_m\supset \mathcal F_{m-1}\supset \dots
\supset \mathcal F_1\supset \mathcal F_0 =0,\]
 such that each factor
$\mathcal E_i$ $:=$ $\mathcal F_i/\mathcal F_{i-1}$ is stable under
$Z_{s_0,t_0}$. For any point $(s,t)$ on $W$, we have the slope
$\mu_{s,t}(\mathcal E_i)$ $=$ $\mu_{s,t}(\mathcal F)$. Otherwise,
$\mathcal F$ is always unstable under $Z_{s_0,t_0\pm \epsilon}$.
Since the actual walls on the second quadrant are nested
semicircular wall, this contradicts the fact that $\mathcal F$ is
locally stable under $Z_{s_0,t_0}$.
\begin{lemma}
Let $\mathcal E_1,\dots,\mathcal E_m$ be the stable factors of
$\mathcal F$ as above, then we have:
\[\text{Hom}(\mathcal E_i,\mathcal E_j[2]) = 0,\]  for all $1\leq
i,j \leq m$. \label{lemma: ext2 vanishing for stable factors}
\end{lemma}
\begin{proof}
In order to apply the same trick as that in Lemma \ref{smoothness
lemma}, we show that each $\mathcal E_i$ is stable on the whole
$W_{(s_0,t_0)}$. First, we show that $\mathcal E_i$ is always in
$\mathcal A_s$. If not so, either  $\mu^{GM}_+($H$^{-1}(\mathcal
E_k))$ or $\mu^{GM}_-($H$^{0}(\mathcal E_k))$ falls into the open
region $W_s$ $:=$ $\{s|(s,t)\in W_{(s_0,t_0)}$ for some $t\}$. If
there exists $\mathcal E_k$ such that $\mu^{GM}_+($H$^{-1}(\mathcal
E_k)$) $=$ $\mu^{GM}(H^{-1}(\mathcal E_k)_{\text{max}})$ $\in$
$W_s$, we may assume
\begin{center}
$s_-$ $:=$ $\mu_+(H^{-1}(\mathcal E_k))$ $=$ max$_{1\leq i\leq m}$
$\{ \mu^{GM}_+($H$^{-1}(\mathcal E_i)) \}$;\\
 $k$ $=$
max$\{i|\mu^{GM}_+($H$^{-1}(\mathcal E_i))$ $=$ $s_-\}$.
\end{center}
 Then when $s$ tends to $s_-$ from the right
along $W_{(s_0,t_0)}$, $\mu_{s,t}$ (H$^{-1}(\mathcal E_k)_{max}$[1])
will go to $+\infty$. Let the quotient of H$^{-1}(\mathcal
E_k)_{max}$[1] $\rightarrow$ $\mathcal E_k$ in $\mathcal A_s$ (for
all $s_{-}<s\leq s_0$) be $\mathcal E'_k$, then there is a map from
$\mathcal F_k$ to $\mathcal E'_k$: $\mathcal F_k\rightarrow \mathcal
E_k\rightarrow \mathcal E'_k$.  By the maximum assumption on
$\mu_+(H^{-1}(\mathcal E_k))$ and $k$, this is surjective for all
$s_{-}<s\leq s_0$. Let the kernel of $\mathcal F_k\rightarrow
\mathcal E'_k$ be $\mathcal F'_k$. Since $\mu_{s,t}(\mathcal
F_{k-1})$ $=$ $\mu_{s,t} (\mathcal F)$, which is bounded on
$W_{(s_0,t_0)}$, and $\mu_{s,t}$ (H$^{-1}(\mathcal F_+)_{max}$[1])
tends to $+\infty$ as $s$ tends to $s_-$, we have
$\mu_{s,t}(\mathcal F'_k)$ $>$ $\mu_{s,t}(\mathcal F)$ when $s$
tends to $s_-$. Since the actual walls of $(1,0,1-n)$ are nested,
$\mathcal F$ must be locally stable along $W_{(s_0,t_0)}$, which
contradicts the inequality $\mu_{s,t}(\mathcal F'_k)$ $>$
$\mu_{s,t}(\mathcal F)$. In a similar way, we get a contradiction
for the case
$\mu_-^{GM}(H^{0}(\mathcal E_k))$ $\in$ $W_s$.\\

Next, we show the stableness of  $\mathcal E_k$. Suppose $k$ is the
maximum number such that $\mathcal E_k$ is not stable for some
$(s,t)$ on $W_{(s_0,t_0)}$. There must be a subobject $\mathcal
E''_k$ of $\mathcal E_k$ in $\mathcal A_s$ for some $(s,t) \in$
$W_{(s_0,t_0)}$ such that $\mu_{s,t}(\mathcal E''_k)$ $>$
$\mu_{s,t}(\mathcal E_k)$. Again let the quotient be $\mathcal E'_k$
and consider the kernel $\mathcal F'_k$ of $\mathcal F$
$\rightarrow$ $\mathcal E'_k$. Then $\mathcal F'$ is a subobject of
$\mathcal F$ and $\mu_{s,t}(\mathcal F'_k)$ $>$ $\mu_{s,t}(\mathcal
F)$, which is a contradiction.\\

Now we repeat the same argument in Lemma \ref{smoothness lemma} for
$\mathcal E_j$ and $\mathcal E_i$. When the $W_{(s_0,t_0)}$ has
radius greater than $3$, since on $W_{(s_0,t_0)}$,
$\mu_{s,t}(\mathcal E_i(-3))$ $=$ $\mu_{s,t}((\mathcal F(-3))$ $<$
$\mu_{s,t}\mathcal F$ $=$ $\mu_{s,t}(\mathcal E_j)$. When the radius
is not greater than $3/2$, $\mathcal E_j$ and $\mathcal E_i(-3)[1]$
are both in $\mathcal A(-k-3)$ (since $\mathcal F(-3)[1]$ is in
$\mathcal A(-k-3)$ and $\mathcal E_i(-3)$ has the same slope of
$\mathcal F(-3)$ along the wall). In either case, Hom($\mathcal
E_j,\mathcal E_i(-3))$ $=$ $0$. We get
\begin{center}Hom $(\mathcal E_i$, $\mathcal E_j $[2]) $\simeq$ Hom
($\mathcal E_j ,\mathcal E_i(-3)$)$^*$ $=$ $0$.\end{center}
\end{proof}
Based on the previous two lemmas, we will study the phenomenon of
wall-crossing in a quiver region. Write $\mathcal F[1]$ as an object
\textbf K in a quiver region $\mathcal A(k)$, for some
$-\sqrt{2n}<k<0$. Let $\overrightarrow {\rho}$ be the character
corresponding to $W_{(s_0,t_0)}$. Then the stable factor filtration
of \textbf K at $\overrightarrow {\rho}$ in $\mathcal A(k)$ is
written as:
\[\textbf K = \tilde{\textbf K}_m
\supset\tilde{\textbf K}_{m-1}\supset \dots\supset \tilde{\textbf
K}_1 \supset \tilde{\textbf K}_0 =0,\]
 where $\tilde{\textbf K}_i$
is $\mathcal F_i[1]$, with \textbf K$_i$ $=$ $\tilde{\textbf
K}_i$/$\tilde{\textbf K}_{i-1}$ = $\mathcal E_i[1]$. Let
$\overrightarrow{\rho}_\pm$ be the character at ($s_0\pm \epsilon$,
$t_0)$. The point in $\mathfrak M^{\overrightarrow
{\rho},ss}(\overrightarrow {n})$ that stands for the S-equivalent
class with stable factors $\{$\textbf K$_1,$ $\dots$, \textbf
K$_m\}$ `blows-up' to two different varieties
V$_{\overrightarrow{\rho}_+}$(\textbf K$_1,$ $\dots$, \textbf K$_m$)
and V$_{\overrightarrow{\rho}_-}$(\textbf K$_1,$ $\dots$, \textbf
K$_m$) in $\mathfrak M^{\overrightarrow {\rho}_+,ss}(\overrightarrow
{n})$ and $\mathfrak M^{\overrightarrow {\rho}_-,ss}(\overrightarrow
{n})$ respectively. Let $S_+$, $S_-$, be the sets defined as follows
respectively:
\begin{center}
$\{\textbf L\in \mathcal A(k)| \textbf L$ has a filtration with
stable factors as a strict subset of $\{$\textbf K$_1,$ $\dots$,
\textbf K$_m\}$ (counting multiples); \textbf L[-1] is  in $\mathcal
A_{s_0+\epsilon}$, stable under $Z_{s_0+ \epsilon, t_0}$ and
$\overrightarrow {l}\cdot \overrightarrow {\rho}_+
>0$ (respectively stable under $Z_{s_0- \epsilon, t_0}$ and $\overrightarrow {l}\cdot \overrightarrow {\rho}_-
>0)\}$,
\end{center}
where $\overrightarrow l$ is the type of \textbf L.
\begin{lemma}
For each \textbf K with type
$(n-\frac{k(k-1)}{2},2n-k^2+1,n-\frac{k(k+1)}{2})$ that is stable
with respect to the character $\rho_-$ and has stable factor
filtration with factors \textbf K$_1,$ $\dots$, \textbf K$_m$ (i.e.
the point \textbf K is in V$_{\overrightarrow{\rho}_-}$(\textbf
K$_1,$ $\dots$, \textbf K$_m$)), one may write it as an extension of
two semistable objects in $\mathcal A(k)$:
\[0\rightarrow \textbf K_+ \rightarrow \textbf K\rightarrow \textbf K_-\rightarrow 0\]
with properties:\\
1. Hom$(\textbf K_+, \textbf L_-)$ $=$ Hom$(\textbf L_-, \textbf
K_+)$ $=$ $0$, for any \textbf
L$_+$ $\in$ $S_+$;\\
2.  Hom$(\textbf K_-, \textbf L_+)$ $=$ Hom$(\textbf L_+, \textbf
K_-)$ $=$ $0$, for any \textbf L$_-$ $\in$ $S_-$. \label{lemma:
properties of extended factors}
\end{lemma}
\begin{proof}
Since  \textbf K[-1] and \textbf L$_-$[-1] are stable under $Z_{s_0-
\epsilon, t_0+}$, and $\overrightarrow {l}\cdot \overrightarrow
{\rho}_+$ $=$ $\overrightarrow {k }\cdot \overrightarrow {\rho}_+$
$>0$, Hom $(\textbf L_-,$ \textbf K) $=$ $0$ for any \textbf L$_-$
$\in$ $S_-$. For any sub-object \textbf K$_+$, we have Hom $(\textbf
L_-,$ \textbf K$_+$) $=$ $0$. Similarly we have Hom$(\textbf K_-,
\textbf L_+)$ $=$ $0$ for any
\textbf L$_+$ $\in$ $S_+$.\\

Start from an extension pair $(\textbf K^0_+,\textbf K^0_-)$ $=$
$(\textbf K_1,\textbf K/\textbf K_1)$ for \textbf K.  If Hom(\textbf
L$_+,$ \textbf K$^0_-)$ $\neq$ $0$ for some \textbf L$_+\in S_+$
(the image is a subobject in  \textbf K$^0_-$), then we make an
adjustment for the pair by moving a subobject \textbf L$^0_+$ in
\textbf K$^0_-$ to \textbf K$^0_+$. Denote this extension pair by
$(\textbf K^{00}_+,\textbf K^{00}_-)$. Then we move a quotient
object \textbf L$^0_-$ in \textbf K$^{00}_+$ to \textbf K$^{00}_-$
if there is any. Denote the new pair by $(\textbf K^1_+,\textbf
K^1_-)$. We may repeat this procedure and get pairs $(\textbf
K^i_+,\textbf K^i_-)_{i\geq 0}$. Denote the type of \textbf K$^i_+$
by $\overrightarrow k^i_+$, it is not hard to see that
$\overrightarrow k^i_+\cdot \overrightarrow {\rho}_+$ is
non-decreasing when $i$ increases. $\overrightarrow k^i_+\cdot
\overrightarrow {\rho}_+$ stop increasing when there is no
adjustment at this step, i.e. $(\textbf K^i_+,\textbf K^i_-)$
satisfies the requirements. Since there are only finite
possibilities for the value of $\overrightarrow k^i_+\cdot
\overrightarrow {\rho}_+$, we always get the extension pair.
\end{proof}
It is immediate that Hom$(\textbf K_\pm,\textbf K_\mp)$ $=$ $0$. Let
$V(\textbf K_+, \textbf K_-)$ be the sub-variety in
V$_{\overrightarrow{\rho}_-}$(\textbf K$_1,$ $\dots$, \textbf K$_m$)
consisting of
 objects that can be written as the extension $0\rightarrow \textbf K_+\rightarrow
 \textbf K \rightarrow \textbf K_-\rightarrow 0$.
\begin{lemma}
 Adopting the notation as above, the
 dimension of $V(\textbf K_+, \textbf K_-)$ is at most
 dimExt$^1(\textbf K_-,\textbf K_+)-$ (dimAut$(\textbf K_+)$ $+$ dimAut$(\textbf K_-)$
$-$ $1$). \label{lemma: boundary of dimension for extended type}
\end{lemma}
\begin{proof}
The extension is given by an element in Hom$(\textbf K_-,\textbf
K_+[1])$. As Hom$(\textbf K_\pm,\textbf K_\mp)$ $=$ $0$, the two
extended objects are isomorphic if they are on the same orbit of the
Aut$(\textbf K_+)$ $\times$ Aut$(\textbf K_-)$ action. To proof the
lemma, we only need to show  that if $f$ $\in$ Hom$(\textbf
K_-,\textbf K_+[1])$ induces a complex in $V(\textbf K_+, \textbf
K_-)$, then the stabilizers of $f$ in Aut$(\textbf K_+)$ $\times$
Aut$(\textbf K_-)$ are the scalars.\\

Let $(g^+,g^-)$ be a stabilizer of $f$. Write \textbf K$_\pm$ as
$\mathcal O(-k-1)\otimes H^\pm_{-1}\xrightarrow{I^\pm} \mathcal
O(-k)\otimes H^\pm_0\xrightarrow{J^\pm} \mathcal O(-k+1)\otimes
H^\pm_1$, then we may represent $g^\pm$ as
$(g^\pm_{-1},g^\pm_{0},g^\pm_{1})$ $\in$ ker$d^0$ $\subset$
Hom$^0(\textbf K^\pm,\textbf K^\pm)$, where $g^\pm_{i}$ $\in$
GL$(H^\pm_i)$. $f$ can be written as $(f_{-1},f_0)$ $\in$
ker$d^1$/im$d^0$ $\subset$ Hom$^1(\textbf K^-,\textbf K^+)$/im$d^0$,
when $f_{-1}\in $ Hom$(\mathcal O(-k-1)\otimes H_{-1}^-,$ $\mathcal
O(-k)\otimes  H_0^+)$ and  $f_{0}\in $ Hom$(\mathcal O(-k)\otimes
H_{0}^-,$ $\mathcal O(-k+1)\otimes  H_1^+)$. Then \textbf K is
written as:
\[
\mathcal O(-k-1)\otimes (H^+_{-1}\oplus
H^-_{-1})\xrightarrow{\begin{bmatrix}I^+ & 0\\ f_{-1} &
I^-\end{bmatrix}} \mathcal O(-k)\otimes (H^+_{0}\oplus
H^-_{0})\xrightarrow{\begin{bmatrix}J^+ & 0 \\ f_{0} &
J^-\end{bmatrix}} \mathcal O(-k+1)\otimes (H^+_{1}\oplus H^-_{1}).
\]
As $(g^+,g^-)$ is a stabilizer, $g^+\circ f\circ (g^-)^{-1} -f$ is
an exact cycle in Hom$^1 (\textbf K^-,\textbf K^+)$ which can be
written as $d^0(s)$ for some $s$ $\in$ Hom$^0 (\textbf K^-,\textbf
K^+)$. Write $\tilde{f}$ for $g^+\circ f\circ (g^-)^{-1}$, we have:
\begin{align*}
& \begin{bmatrix} \text{Id}H^+_0 & s_0\\ 0 & \text{Id}H^-_0
\end{bmatrix}
\begin{bmatrix} g^+_0 & 0\\ 0 &  g^-_0
\end{bmatrix}
\begin{bmatrix}
I^+ & f_{-1}\\ 0 & I^- \end{bmatrix}
 =
 \begin{bmatrix} \text{Id}H^+_0 & s_0\\ 0 & \text{Id}H^-_0
\end{bmatrix}
\begin{bmatrix}
g^+_0 I^+ & g^+_0 f_{-1}\\ 0 & g^-_{0}I^- \end{bmatrix}  \\
 = &
\begin{bmatrix} \text{Id}H^+_0 & s_0\\ 0 & \text{Id}H^-_0\end{bmatrix}
\begin{bmatrix} I^+g^+_{-1} &
 \tilde{f}_{-1} g^-_{-1}\\ 0 & I^-g^-_{-1} \end{bmatrix}
  =
\begin{bmatrix} I^+g^+_{-1} &
 (\tilde{f}_{-1} + s_0I^-) g^-_{-1}\\ 0 & I^-g^-_{-1}
 \end{bmatrix}\\
  = &
\begin{bmatrix} I^+g^+_{-1} &
 (f_{-1} + I^-s_{-1}) g^-_{-1}\\ 0 & I^-g^-_{-1} \end{bmatrix}
 =
 \begin{bmatrix}
I^+ & f_{-1}\\ 0 & I^- \end{bmatrix}
\begin{bmatrix} \text{Id}H^+_{-1} & s_{-1}\\ 0 & \text{Id}H^-_{-1}
\end{bmatrix}
\begin{bmatrix}
g^+_{-1} & 0\\ 0 & g^-_{-1}\end{bmatrix}.
\end{align*}
By changing the labels, we have a similar equality for $J$,
$\begin{bmatrix} \text{Id}H^+ & s\\ 0 & \text{Id}H^-
\end{bmatrix}
\begin{bmatrix}
g^+ & 0\\ 0 & g^-\end{bmatrix}$ is now a morphism from \textbf K to
itself. Since \textbf K is stable, Hom(\textbf K,\textbf K)$=$ $\mathbb C$,
$\begin{bmatrix} g^+ & 0\\
0 & g^-\end{bmatrix}$ must be identity up to a scalar.
\end{proof}
Let the type of \textbf K$_+$ be $\overrightarrow {n}_+$, we write
$\mathfrak M^{\overrightarrow{\rho_0},ss}(\overrightarrow {n}_+)$ as
the moduli space of semistable objects in $\mathcal A(k)$ with type
$\overrightarrow {n}_+$ and character $\overrightarrow{\rho_0}$. It
is a projective variety as the case of Hilbert case. For a positive
integer $c$, let $\mathfrak
M_c^{\overrightarrow{\rho_0},ss}(\overrightarrow {n}_+)$ be the
locus in $\mathfrak M^{\overrightarrow{\rho_0},ss}(\overrightarrow
{n}_+)$ consisting of points that dimHom(\textbf K$_+$, \textbf
K$_+$) $=$ $c$. As the constrain is algebraic, $\mathfrak
M^{\overrightarrow{\rho_0},ss}(\overrightarrow {n}_+)$ is a
subvariety in $\mathfrak
M^{\overrightarrow{\rho_0},ss}(\overrightarrow {n}_+)$, and we may
stratify $\mathfrak M^{\overrightarrow{\rho_0},ss}(\overrightarrow
{n}_+)$ as $\sqcup_{c\in\mathbb N} \mathfrak
M_c^{\overrightarrow{\rho_0},ss}(\overrightarrow {n}_+)$. These
notations also make sense for $\overrightarrow {n}_-$
respectively.\\

Let $V_{-,\overrightarrow {\rho_0}}(\overrightarrow
{n}_+,\overrightarrow {n}_-)$ $\subset$ $\mathfrak
M^{\overrightarrow {\rho}_-,s}(\overrightarrow {n})$ be the locus
consisted of object that can be extended by \textbf K$_+$ and
\textbf K$_-$ which satisfy the properties in Lemma \ref{lemma:
properties of extended factors} with type $\overrightarrow {n}_+$
and $\overrightarrow {n}_-$. Let $V_{-,\overrightarrow {\rho_0}}$
$\subset $ $\mathfrak M^{\overrightarrow {\rho}_-,s}(\overrightarrow
{n})$ consisting of objects that are stable with respect to
$\overrightarrow {\rho}_-$ but not $\overrightarrow {\rho}_+$, then
by Lemma \ref{lemma: properties of extended factors},
$V_{-,\overrightarrow {\rho_0}}$ $=$ $\bigcup_{\overrightarrow {n}_+
+ \overrightarrow {n}_- = \overrightarrow {n}} V_{-,\overrightarrow
{\rho_0}}(\overrightarrow {n}_+,\overrightarrow {n}_-)$, we may
estimate the dimension of $V_{-,\overrightarrow {\rho_0}}$ by
studying the dimension of each piece.
\begin{proposition}
Adopting the notation as above, when $-\sqrt{2n}<s<0$, the dimension
of $V_{-,\overrightarrow {\rho_0}}$ is less than $2n-2$.
\label{contracting is more than producing lemma}
\end{proposition}

\begin{proof}
For $c,d$ $\in$ $\mathbb N$, let $V_{-,\overrightarrow
{\rho_0}}(\overrightarrow {n}_+,\overrightarrow {n}_-)(c,d)$ be the
locus where the complex can be extended by objects in
$\tilde{\mathfrak M}_c^{\overrightarrow{\rho_0},ss}(\overrightarrow
{n}_+)$ and $\tilde{\mathfrak
M}_d^{\overrightarrow{\rho_0},ss}(\overrightarrow {n}_-)$. Then we
have
 \begin{align*}
&\text{dim}V_{-,\overrightarrow {\rho_0}}(\overrightarrow
{n}_+,\overrightarrow {n}_-)(c,d) \\
 \leq & \text{
dim}\tilde{\mathfrak
M}_c^{\overrightarrow{\rho_0},ss}(\overrightarrow {n}_+) + \text{
dim}\tilde{\mathfrak
M}_d^{\overrightarrow{\rho_0},ss}(\overrightarrow {n}_-) +
\text{max}_{\textbf K_+,\textbf K_-}\{\text{ dimExt}^1(\textbf
K_-,\textbf K_+)\} - c - d +1 \tag {by Lemma \ref{lemma: boundary of
dimension for
extended type}}\\
\leq & \text{dimExt}^1(\textbf K_+,\textbf K_+) +
\text{dimExt}^1(\textbf K_-,\textbf K_-) + \text{dimExt}^1(\textbf
K_-,\textbf K_+) - c -d  +1 \tag {the dimension of Zariski tangent space}\\
 = & -\chi(\textbf K_+,\textbf K_+)-\chi(\textbf K_-,\textbf
K_-)-\chi(\textbf K_-,\textbf K_+) +1 \tag {by Lemma \ref{lemma: ext2 vanishing for stable factors} and \ref{lemma: properties of extended factors}}\\
 = & -\chi(\textbf K,\textbf K) + \chi(\textbf K_+,\textbf K_-) + 1\\
 = & 2n + \chi(\textbf K_+,\textbf K_-).
\end{align*}

The remaining task is to estimate $\chi(\textbf K_+,\textbf K_-)$.
Since when $s$ moves from $s_0+\epsilon$ to $s_0-\epsilon$,
$\overrightarrow {\rho}$ will move from $\overrightarrow {\rho}_0$
$+$ $\epsilon(n_{(0)},-n_{(-1)},0)$ to $\overrightarrow {\rho}_0$
$-$ $\epsilon(n_{(0)},-n_{(-1)},0)$, and \textbf K$_-$ does not
destabilize $\textbf K$ on the $s_0+\epsilon$ side, we have
\begin{align*}
\overrightarrow {n}_+\cdot (n_{(0)},-n_{(-1)},0) > 0,
\overrightarrow {n}_+\cdot (0,n_{(1)},-n_{(0)}) >0.
\tag{*}\label{inequality: n+ and rho+}
\end{align*}
Now we have the estimation on $\chi (\textbf K_+,\textbf K_-)$:
\begin{align*}
&  -\chi (\textbf K_+,\textbf K_-) \geq \chi (\textbf K_-,\textbf K_+) - \chi (\textbf K_+,\textbf K_-)\\
= & (\overrightarrow{n}_-\cdot \overrightarrow{n}_+
-3(n_{-(-1)}n_{+(0)}+n_{-(0)}n_{+(1)}) +6n_{-(-1)}n_{+(1)} )-
(\overrightarrow{n}_-\cdot \overrightarrow{n}_+
-3(n_{+(-1)}n_{-(0)} \\ &  +n_{+(0)}n_{-(1)}) +6n_{+(-1)}n_{-(1)})\\
= & (6n_{(-1)}-3n_{(0)})n_{+(1)} +(3n_{(1)}-3n_{(-1)})n_{+(0)}
+(3n_{(0)}-6n_{(1)})n_{+(-1)} \\
= & 3\overrightarrow{n}_+\cdot (k-1,-k,k+1)\\
\geq & 3.
\end{align*}
The last inequality is due to $(k-1,-k,k+1)$ $=$
$\frac{k-1}{n_{(0)}}(n_{(0)},-n_{(-1)},0)$ $-$
$\frac{k+1}{n_{(0)}}(0,n_{(1)},-n_{(0)})$ and the formula
(\ref{inequality: n+ and rho+}).
\end{proof}
\begin{corollary}
For a generic $\overrightarrow {\rho}$ not on any actual
destabilizing wall, $\mathfrak M^{\overrightarrow
{\rho},s}(\overrightarrow {n})$ $=$ $X$//$_{\rho}G$ is irreducible
and smooth. \label{smoothness and irred property}
\end{corollary}
\begin{proof}
In the Hilbert scheme chamber, the irreducible components that
contain $X^{\overrightarrow {\rho},s}$ are reduced and irreducible
since Hilb$^nS$ is so. Passing to another quiver region does not
affect the reduceness property of $X^{\overrightarrow {\rho},s}$. By
Proposition \ref{contracting is more than producing lemma}, while
going across one destabilizing wall, the new stable locus
$V_{-,\rho_0}$ in $X^{\overrightarrow{\rho}_-,s}$//$G$ has dimension
less than $2n-2$. Therefore the dimension of
$X^{\overrightarrow{\rho}_-,s}\setminus
X^{\overrightarrow{\rho}_+,s}$ is less than $2n-2+$ dim$G$. On the
other hand, the total space $X$ is  Spec$\mathbb C[M]/(J\circ I)$,
where $M$ is the space Hom$(H_{-1},H_0)$ $\otimes$ Hom$(\mathcal
O(-k-1),\mathcal O(-k))$ $\times$ Hom $(H_{0},H_1)$ $\otimes$
Hom$(\mathcal O(-k),\mathcal O(-k+1))$. The dimension of $M$ is
$3n_{-1}n_0+3n_0n_1$, and $J\circ I$ has $6n_{-1}n_1$ equations. In
any quiver region, we have $3n_{-1}n_0+3n_0n_1-6n_{-1}n_1$ $=$
$2n-2$ $+$ dim$G$. Each irreducible component has dimension at least
$2n+$ dim$G$. Since $X^{\overrightarrow{\rho}_-,s}$ is open in $X$,
and dim$X^{\overrightarrow{\rho}_-,s}\setminus
X^{\overrightarrow{\rho}_+,s}$ $<$ $2n-2+$ dim$G$, we get
$X^{\overrightarrow{\rho}_-,s}$ is irreducible.\\

The dimension of the Zariski
tangent space at a point $\textbf K =(I_0,J_0)$ is the dimension of
 Hom$_{\mathbb C}$ ($\mathbb C[M]/(J\circ I), \mathbb C[t]/(t^2))$ at
 point $(I_0,J_0)$. Each tangent
 direction is written in a form $(I_0,J_0)$ $+$ $t(I_1,J_1)$. In
 order to
 satisfy the equation $J\circ I =0$, we need
 \[J_0\circ I_1+J_1\circ  I_0 =0.\]
Hence the space of $(I_1,J_1)$ is just the kernel of $d^1:$
Hom$^1(\textbf K,\textbf K)$ $\rightarrow$ Hom$^2(\textbf K,\textbf
K)$. Now by Lemma \ref{smoothness lemma}, $d^1$ is surjective. The
Zariski tangent space has dimension dimHom$^1(\textbf K,\textbf K)$
$-$ dimHom$^2(\textbf K,\textbf K)$ $=$
$3n_0(n_1+n_{-1})-6n_{-1}n_1$, which is the dimension of $M$ minus
the number of equations. We get the smoothness of $X^{s,\rho}$.
Furthermore, since $(n_{-1},n_0,n_{1})$ $=$ $(n-\frac{(k-1)k}{2}, 2n
-k^2 +1, n-\frac{(k+1)k}{2})$ is primitive, $G$ acts freely on
$X^{s,\rho}$. By Luna's \'{e}tale slice theorem, $X^{s,\rho}$
$\rightarrow$ $ X$//$_{\rho}G$ is a principal bundle. Since
$X^{s,\rho}$ is smooth, by Proposition IV.17.7.7 in \cite{Gr}, the
base space is also smooth.
\end{proof}
\subsection{properties of GIT}
Birational geometry via GIT has been studied in \cite{DoHu} by
Dolgachev and Hu, \cite{Th} by Thaddeus. In this section, we
recollect some properties in a
language of the affine GIT.\\

Let $X$ be an affine algebraic $G$-variety , where $G$ is a
reductive group and acts on $X$ via a linear representation. Given a
character $\rho$: $G\rightarrow \mathbb C^{\times}$, the
(semi)stable locus is written as $X^{s,\rho}$ ($X^{ss,\rho}$). We
write $C[B]^{G,\chi}$ for the $\chi$-semi-invariant functions on
$B$, i.e. one has
\[f(g^{-1}(x))=\chi(g)\cdot f(x),\text{ } \forall g\in G, x\in B.\]
Denote the GIT quotient by $X$//$_\rho G$ $:=$ \textbf{Proj}
$\bigoplus_{n\geq 0} \mathbb C[X]^{G,\rho^n}$ and the map from
$X^{ss,\rho}$ to $X$//$_\rho G$ by $F_{\rho}$.\\

In additions, we need he following assumptions on $X$ and $G$: 1.
there are only finite many walls in the space of characters on which
there are semistable but non-stable points, in the chamber we have
$X^{s,\rho}$ $=$ $X^{ss,\rho}$. 2. $X^{s,\rho}$ is smooth and the
action of $G$ on $X^{s,\rho}$ is free. 3. $\mathbb C[X]^G$ $=$
$\mathbb C$, i.e., $X$//$_\rho G$ is projective and connected. 4.
The closure of $X^{s,\rho}$ (if non-empty) for any $\rho$ is the
same irreducible component. 5. Given any point $x$ $\in$ $X$, the
set of characters $\{\rho|$ $x$ $\in$ $X^{ss,\rho}\}$ is closed.\\

Notations and constructions: let $\rho$ be a generic character (i.e.
not on any walls) satisfying that $X^{s,\rho}$ is non-empty, then by
assumptions we have a $G$-principal bundle $X^{s,\rho}$ $\rightarrow
$ $X$//$_\rho G$ $=$ $X^{s,\rho}$/$G$. Giving another character
$\rho_0$ of $G$, we denote $\mathcal L_{\rho,\rho_0}$ to be the line
bundle over $X$//$_\rho G$ by composing the transition functions of
the $G$-principal bundles with $\rho_0$. Now we are ready to list
some properties from the variation geometric invariant theory
(VGIT).
\begin{proposition}
Let $X$ be an affine algebraic $G$-variety that satisfies the
assumptions 1 to 5 and  $\mathcal L_{\rho,\rho_0}$ be as defined
above. We
have:\\
1. $\Gamma$ $(X$//$_\rho G$, $\mathcal
L_{\rho,\rho_0}^{\otimes n})$ $\simeq$ $\mathbb
C[X^{s,\rho}]^{G,\rho_0^n}$.\\
2. If $\rho_+$ and $\rho$ are in the same chamber, then $\mathbb
C[X^{s,\rho}]^{G,\rho_+^n}$ $=$ $\mathbb C[X]^{G,\rho_+^n}$ for
$n\gg 1$, $\mathcal L_{\rho,\rho_+}$ is ample; if $\rho_0$ is a
generic point on the wall of the $\rho$-chamber,
then $\mathcal L_{\rho,\rho_0}$ is nef and semi-ample.\\
3. Let $\rho_+$ and $\rho_0$ be in the chamber of $\rho$ and on the
wall respectively, then there is an inclusion $X^{ss,\rho_+}$
$\subset$ $X^{ss,\rho_0}$ inducing a canonical projective morphism
pr$_+$: $X$//$_{\rho_+}G$
$\rightarrow$ $X$//$_{\rho_0}G$.\\
4. A curve $C$ (projective, smooth, connected) in $X$//$_{\rho_+}G$
is contracted by pr$_+$ if and only if it is contracted by
$X$//$_{\rho_+}G$ $\rightarrow$ \textbf{Proj}$\oplus_{n\geq0}
\Gamma(
X$//$_{\rho_+}G,\mathcal L_{\rho_+,\rho_0}^{\otimes n})$.\\
5. Let $\rho_+$ and $\rho_-$ be in two chambers on different sides
of the wall, let $\rho_0$ be a generic point on the wall. Assume
that $X^{s,\rho_{\pm}}$ are both non-empty, then the morphisms
$X$//$_{\rho_{\pm}}G$ $\rightarrow$ $X$//$_{\rho_0}G$ are proper and
birational. If they are both small, then the rational map
$X$//$_{\rho_{-}}G$ $\dashrightarrow$ $X$//$_{\rho_{+}}G$ is a flip
with respect to $\mathcal L_{\rho_+,\rho_0}$. \label{recollections
from VGIT}
\end{proposition}
\begin{proof} 1. This is true for general $G$-principal bundle by flat
descent theorem, see \cite{SGA} Expos\'{e} I, Th\'{e}or\`{e}me 4.5.\\

2 and 3. By assumption 5,  $X^{s,\rho}$ $\subset$ $X^{ss,\rho_*}$.
By assumption 4, the natural maps: $\mathbb C[X]^{G,\rho_*^n}$
$\rightarrow$ $\mathbb C[X^{s,\rho}]^{G,\rho_*^n}$ $\simeq$ $\Gamma$
$(X$//$_\rho G$, $\mathcal L_{\rho,\rho_*}^{\otimes n})$ in
injective for $* =0,+$ and $n\in \mathbb Z_{\geq 0}$. Hence the base
locus of $\mathcal L_{\rho,\rho_*}$ is empty. $\mathcal
R(X$//$_{\rho}G$, $\mathcal L_{\rho,\rho_*})$ $\simeq$
$\bigoplus_{n\geq 0} \mathbb C[X^{s,\rho}]^{G,\rho_*^n}$ is finitely
generated over $\mathbb C$. The canonical morphism $X$//$_{\rho}G$
$\rightarrow$ \textbf{Proj}$\bigoplus_{n\geq 0} \mathbb
C[X^{s,\rho}]^{G,\rho_*^n}$ is birational and projective when
$X^{s,\rho_*}$ is non-empty. Now we have series of morphisms:
\begin{center}
pr$_+$: $X$//$_{\rho}G$ $\rightarrow$ \textbf{Proj}$\bigoplus_{n\geq
0} \mathbb C[X^{s,\rho}]^{G,\rho_*^n}$ $\rightarrow$
\textbf{Proj}$\bigoplus_{n\geq 0} \mathbb C[X]^{G,\rho_*^n}$ $=$
$X$//$_{\rho_*}G$.
\end{center}
The morphism pr$_+$ maps each $\rho_*$ S-equivariant class to itself
set-theoretically. When $\rho_+$ is in the same chamber of $\rho$,
by the assumption 2, this is an isomorphism, implying that $\mathcal
L_{\rho,\rho_+}$ must be ample and $\mathbb
C[X^{s,\rho}]^{G,\rho_*^n}$ $=$ $\mathbb C[X]^{G,\rho_*^n}$ for $n$
large enough. By the definition of $\mathcal L_{\rho,\rho_+}$, it
linearly extends to a map from the space of $\mathbb R$-characters
of $G$ to NS$_{\mathbb R}$($X$//$_{\rho}G$). Since all elements in
the $\rho$ chamber are mapped
into the ample cone, $\rho_0$ must be nef.\\

4. `$\Leftarrow$' direction:  by the assumption 4,
\textbf{Proj}$\oplus_{n\geq0} \mathbb C[X^{s,\rho_+}]^{G,\rho_0^n}$
$\rightarrow$ \textbf{Proj}$\oplus_{n\geq0} \mathbb
C[X]^{G,\rho_0^n}$ is always surjective. If $C$ is contracted at
\textbf{Proj}$\oplus_{n\geq0} \mathbb C[X^{s,\rho_+}]^{G,\rho_0^n}$,
then it is also contracted at \textbf{Proj}$\oplus_{n\geq0} \mathbb
C[X]^{G,\rho_0^n}$.\\

`$\Rightarrow$': Let $G'$ be the kernel of $\rho_0$, we show that
there is a
subvariety $P$ in $X^{s,\rho_+}$ satisfying:\\
A. $P$ is a $G'$-principal bundle, and the base space is projective, connected;\\
B. $F_{\rho_+} (P)$ $=$ $C$.\\
Suppose we find such $P$, then any function $f$ in $\mathbb
C[X^{s,\rho_+}]^{G,\rho_0^n}$ is constant on each $G'$ fiber. Since
the base space is projective and connected, it must be constant on
the whole space $P$. Since $F_{\rho_+} (P)$ $=$ $C$, the value of
$f$ on $F_{\rho_+}^{-1}(C)$ is determine by this constant. Hence the
canonical morphism contracts $C$ to a point.\\

To get $P$, we may assume $G'$ $\neq$ $G$, choose $N$ large enough
and finitely many $f_i$'s in $\mathbb C[X]^{G,\rho_0^N}$ such that
$\bigcap_i (V(f_i) \cap$ $F_{\rho_0}^{-1}(pr_+(C))$ is empty. Since
all points in $F_{\rho_0}^{-1}(pr_+(C))$ are S-equivariant in
$X^{ss,\rho_0}$, each $\overline{Gx}$ contains all minimum orbits
$Gy$ in $F_{\rho_0}^{-1}(pr_+(C))$. Choose $y$ such that $Gy$ is
closed in $X^{ss,\rho_0}$, let $P_y$ be
\[\bigcap_{i}\{x\in F_{\rho_+}^{-1}(C)|f_i(x) =f_i(y)\}.\] For any
$p$ $\in$ $C$, since the $G$-orbit $\overline{F_{\rho_+}^{-1}(p)}$
contains $y$ and $G$ is reductive, there is a subgroup $\beta$:
$\mathbb C^\times$ $\rightarrow$ $G$ and $x_p$ $\in$
$F^{-1}_{\rho_+}(p)$ satisfying that $y$ $\in$
$\overline{\beta(\mathbb C)\times x_p}$. Since $y$ $\in$
$X^{ss,\rho_0}$, there is a $\rho_0^N$-semi-invariant $f_i$ such
that $f_i(y)$ $=$ $0$. Therefore $\rho_0\circ \beta$ $=$ $0$. This
implies that for any $\rho_0$-semi-invariant function $f$ $f(x_p)$
$=$ $f(y)$. Condition b is checked. Let $G''$ be the kernel of
$\rho^N$. By the choices of $f_i$'s, another point $x_q$ on $Gx_p$
is in $P_y$ if and only if they are on the same $G''$-orbit. Since
$G$ acts freely on all stable points, $P_y$ becomes a $G''$
principal bundle over base $C$. As $[G'':G']$ is finite, we may
choose a connected component of $P_y$ and as a $G'$-principal
bundle, the induced morphism from base
space to $C$ is finite. Condition A is checked.\\

5. This is due to Theorem 3.3 in \cite{Th}.
\end{proof}
\begin{remark}
When the difference between $X^{s,\rho_+}$ and $X^{s,\rho_-}$ is of
codimension two in $X^{s,\rho_+}\cup X^{s,\rho_-}$, since
$X^{s,\rho_+}\cup X^{s,\rho_-}$ is smooth, irreducible and
quasi-affine by the second assumption, we have:
\begin{center}
$\mathbb C[X^{s,\rho_+}]^{G,\rho_-^n}$ $=$ $\mathbb
C[X^{s,\rho_+}\cup X^{s,\rho_-}]^{G,\rho_-^n}$ $=$ $\mathbb
C[X^{s,\rho_-}]^{G,\rho_-^n}$ $=$ $\mathbb C[X]^{G,\rho_-^n}$ for
$n\gg 0$.
\end{center}
In this case,  the rational morphism between $X^{s,\rho_+}$ and
$X^{s,\rho_-}$ identifies NS$_{\mathbb R} (X$//$_{\rho_+}G$) and
NS$_{\mathbb R} (X$//$_{\rho_-}G$). It maps $[\mathcal
L_{\rho_+,\rho_*}]$ to $[\mathcal L_{\rho_-,\rho_*}]$ for all
$\rho_*$ in either $\rho_\pm$ chamber. \label{remark on the divisor
glueing in flip case}
\end{remark}
\subsection{Walls on the Second Quadrant}
Now the correspondence picture of the stable base locus
decomposition of the effective cone and the actual destabilizing
walls in the second quadrant is clear:
\begin{theorem}
In the second quadrant of the $(s,t)$-plane of Bridgeland stability
conditions, the semicircular actual walls in
$W^{\text{actual}}_{(1,0,1-n)}$ is one to one corresponding to
stable base locus decomposition walls on one side of the divisor
cone of Hilb$^nS$. \label{left half upper plane's main theorem in
the body}
\end{theorem}
\begin{proof}
Each point in $\{(t,s)|$ $0<t<\frac{1}{2},-\sqrt{2n}<s<0\}$ falls
into some quiver region $\mathcal A(k)$. As explained before
Proposition \ref{finite many wall, t>0 is HilbS}, the moduli space
of $Z_{s,t}$-semistable objects with invariants $(r,c_1,\chi)$ $=$
$(1,0,1-n)$ is parameterized by the quotient space
$X_k$//$_{\overrightarrow{\rho}_{s,t,k}}G_k$. By Proposition
\ref{finite many wall, t>0 is HilbS}, there are finitely many actual
destabilizing walls, and in each chamber the moduli space remains
the same.
By the formula \ref{character rho of s,t}, the character
$\overrightarrow {\rho}_{s,t,k}$ $=$ $(\rho_{-1},\rho_0,\rho_1)$
always satisfies $\rho_{-1}>0>\rho_1$.\\

We first check that the $G_k$-variety $X_k$ satisfies the
assumptions of Proposition \ref{recollections from VGIT} for all
$\overrightarrow {\rho}_{s,t,k}$. The assumption 1 `finiteness of
walls' is due to the second property in Proposition \ref{finite many
wall, t>0 is HilbS}. The smoothness and irreducible property is
checked in Corollary \ref{smoothness and irred property}. $\mathbb
C[X_0]^{G_0}$ $=$ $\mathbb C$, since Hilb$^nS$ is projective. Since
$\mathbb C[X_k]^{G_k}$ $=$ $\mathbb C$ if and only if for some
$\overrightarrow {\rho}_{s,t,k}$, $X_k$//$_{\overrightarrow
{\rho}_{s,t,k}}G_k$ is projective or empty, this is checked by
induction on $k$. The last assumption 5
holds by King's  criterion \cite{Ki} for (semi)stable quiver
representation. \\

Now we may assign a divisor $[\mathcal
L_{\overrightarrow{\rho_+},\overrightarrow{\rho}_{s,t,k}}]$ to
$X_k$//$_{\overrightarrow {\rho}_{s,t,k}}G_k$, where
$\overrightarrow {\rho}_+$ is the character in the chamber.
Starting from a sufficient small $t>0$ and $-1<s<0$, at where
$X_0$//$_{\overrightarrow {\rho}_{s,t,0}}G_0$ is Hilb$^nS$, let $t$
fix and $s$decrease. At an actual destabilizing wall, let pr$_+$ be
the morphism from $X_k$//$_{\overrightarrow
{\rho}_{s_0+\epsilon,t,k}}G_k$ to $X_k$//$_{\overrightarrow
{\rho}_{s_0,t,k}}G_k$ as that in Proposition \ref{recollections from
VGIT}.
One of three different cases may happen:\\
1. pr$_+$ is a small contracton;\\
2. pr$_+$ is birational and has an exceptional divisor;\\
3. all points in Hilb$^nS$ are destabilized.\\

Now by Proposition \ref{contracting is more than producing lemma},
in Case 1, we get small morphism on both sides. In addition,
$X_k^{s,\rho_+}\setminus X_k^{s,\rho_-}$ has codimension not less
than $2$, else
$X_k$//$_{\overrightarrow{\rho}_{s_0-\epsilon,t,k}}G_k$
cannot be projective. 
By property 5 in Proposition \ref{recollections from VGIT}, this is
the flip with respect to the divisor $[\mathcal
L_{\overrightarrow{\rho_+},\overrightarrow{\rho}_{s_0,t,k}}]$. As
the different part of  $X_k^{s,\rho_+}$ and  $X_k^{s,\rho_-}$ is of
codimension $2$, their divisor cones are identified as explained in
Remark \ref{remark on the divisor glueing in flip case}. While $s$
decreases 
$\rho_1/\rho_{-1}$ is decreasing, so the divisor always jumps to the
next chamber and does not go back. \\

In Case 2, $X_k$//$_{\overrightarrow {\rho}_{s_0-\epsilon,t,k}}G_k$
$\rightarrow$ $X_k$//$_{\overrightarrow {\rho}_{s_0,t,k}}G_k$ does
not have any exceptional divisor by Proposition \ref{contracting is
more than producing lemma}, hence the Picard number of
$X_k$//$_{\overrightarrow{\rho}_{s_0-\epsilon,t,k}}G_k$  is $1$. By
property 4 in Proposition \ref{recollections from VGIT}, Case 2 only
happens when the the canonical model associate to $\mathcal
L_{\overrightarrow{\rho_+},\overrightarrow{\rho}_{s_0,t,k}}$
contracts a divisor, i.e. the identified divisor of $\mathcal
L_{\overrightarrow {\rho}_+,\overrightarrow {\rho}_{s,t,k}}$ on
Hilb$^nS$ is on the boundary of the Movable cone. The next
destabilizing wall on the left corresponds to the zero divisor, it
must be Case 3. In general, if the boundary of the Movable cone is
not the same as that of the Nef cone, then Case 2 happens.
Otherwise, case 2 does not happen and the procedure ends up with a
Mori fibration of Case 3.\\

Besides all previous ingredients, we only need to check that Case 3
happens before $s=-\sqrt{2n}$ when $t=0+$.
\begin{lemma}
There is a semicircular wall with radius greater than $1$ such that
inside the wall, there is no semistable object with invariant
$(1,0,1-n)$. \label{lemma: existence of the last wall}
\end{lemma}
\begin{proof}
When $(k+2)(k+1)>2n$, $\mathcal O(-k)[1]$ always has non-zero map to
any object $\mathcal A(-k)$ with invariant $(n_{-1},n_0,n_1)$ $=$
$(n-\frac{(k-1)k}{2},2n-k^2+1,n-\frac{(k+1)k}{2})$, since $2n-k^2+1$
$>$ $3(n-\frac{(k+1)k}{2})$. $\mathcal O(-k)[1]$ corresponds to the
potential wall across $(-k,0)$, hence there is no stable object with
invariant $(1,0,1-n)$ inside this semicircle.
\end{proof}
By the lemma, Case 3 must happen on this wall or a larger actual
wall.
\end{proof}
\subsection{The Vertical Wall and the First Quadrant}
\begin{proposition}
Suppose $\sigma$ of the Sklyanin algebra
Skl$(E,\mathcal L,\sigma)$
 is of infinite order, then no curve is contracted on the vertical wall $s=0$, i.e., the
vertical wall is a faked wall.
\end{proposition}
\begin{proof}
By the formula in Example \ref{character rho of s,t for k=0}, the
vertical wall corresponds to the wall in $\mathcal A(0)$ with
respect to the character $(1,0,-1)$ . We need the following
criterion for the stable monad. 
\begin{lemma} Suppose $\sigma$ of the Sklyanin algebra
Skl$(E,\mathcal L,\sigma)$
 is of infinite order, then a monad \textbf K:
$\mathcal O(-1)\otimes \mathbb C^{n}\rightarrow \mathcal O \otimes
\mathbb C^{2n+1}\rightarrow \mathcal O(1) \otimes \mathbb C^{n}$ is
stable with respect to $(1,0,-1)$ if and only if the first map is
injective, the second map is surjective and the homological sheaf
H$^0(\textbf K)$ at the middle term is a line bundle.
\end{lemma}
\begin{proof}[Proof of lemma] By the discussion in
Proposition \ref{finite many wall, t>0 is HilbS}, character
$(1,0,-1)$ is on the wall of the `hilbert scheme' chamber that
contains $((2n+1)(m-1),n,-(2n+1)m)$, for $m\gg 1$, hence any
$(1,0,-1)$-stable point is $((2n+1)(m-1),n,-(2n+1)m)$-stable and
corresponds to a sheaf of invariant $(r,c_1,\chi)$ $=$ $(1,0,1-n)$.
Denote $I$ and $J$ as the map from $\mathcal O(-1)\otimes \mathbb
C^n$ to $\mathcal O \otimes \mathbb C^{2n+1}$ and from $\mathcal O
\otimes \mathbb C^{2n+1}$ to $\mathcal O(1)\otimes \mathbb C^n$
respectively. Write $I$ $=$ $xI_1 +yI_2 +zI_3$, where $I_k$ is a
linear map from $\mathbb C^{n}$ to $\mathbb C^{2n+1}$, then the
monad corresponds to a line bundle if and only if the cokernel of
$I$ is a vector bundle . By Corollary 3.12 and Lemma 3.11 in
\cite{NS} on the criterion of vector bundle, by restricting on $E$,
H$^0(\textbf K)$ is a line bundle if and only if $lI_1 + mI_2 +nI_3$
is injective for all non-zero triple $(l,m,n)$ $\in$ $\mathbb C^3$
(or equivalently all $[(l,m,n)] \in E$). Now we may show the `if'
and `only if' statements.
\\\\
`$\Rightarrow$': Suppose H$^0(\textbf K)$ is not a line bundle, then
$lI_1 + mI_2 +nI_3$ has a non-zero element $v_{-1}$ in its kernel,
then we consider the subcomplex that generated by $v_{-1}$, i.e. the
minimum subcomplex that contains $v_{-1}$. It is not hard to check
that this subcomplex has dim$(H'_{-1},H'_0,H'_1)$ $=$ $(1,0,0)$,
$(1,1,0)$, $(1,2,0)$ or $(1,2,1)$. Either case contradicts the
$(1,0,-1)$-stableness
requirement.\\\\
`$\Leftarrow$': Suppose the complex is not $(1,0,-1)$ stable, then a
subcomplex with type $(a,b,c)$ destabilizes the monad. Since \textbf
K is $((2n+1)(m-1),n,-(2n+1)m)$-stable for $m\gg 1$,  we have $b$
$\leq$ $2a$ $=$ $2c$. Restricting on the elliptic curve $E$, since
$I$ is injective at every point, we have a complex on $E$:
\[0\rightarrow \mathcal L^* \otimes \mathbb C^{a}\rightarrow
\mathcal O_E\otimes \mathbb C^b\rightarrow \overline{\mathcal
L}\otimes \mathbb C^a\rightarrow 0,\] which is exact except the
middle term. Comparing the rank and the degree, we get $b$ $=$ $2a$
and the complex is exact. But that is not possible since $\mathcal
L^{*\otimes a}\otimes \overline{\mathcal L}^{\otimes a}$ $\nsimeq$
$\mathcal O$ given that $\sigma^{3a}$ is not id$_E$.
\end{proof}
Back to the proof: According to the proof of the lemma, any complex
whose H$^0(\textbf K)$ is not a line bundle has $(1,2,1)$-type
subcomplex, hence each $(1,0,-1)$-semistable complex has a
filtration with $(1,0,-1)$-stable factors of the following types:
one copy of $(a,2a+1,a)$ (a line bundle $\mathcal E$) and several
$(1,2,1)$'s (quotient points
$\mathcal O_p[-1]$ for $p$ $\in$ $E$).\\

Basic computation shows that: Ext$^1$ ($\mathcal E$, $\mathcal
O_p[-1]$) is $0$; Ext$^1$ ($\mathcal O_p[-1]$, $\mathcal E$) is
$\mathbb C$ for all $p$; Ext$^1$ ($\mathcal O_p[-1]$, $\mathcal
O_q[-1]$) is $\mathbb C$ if and only if $p=q$ or
$p=\sigma^3(q)$ and is $0$ for any other $q$.\\
Hence dimExt$^1$ of any two factors is at most dimension $1$, and
any S-equivariant class has only finitely many non-isomorphic
complexes, which means no curve is contracted.
\end{proof}
\begin{lemma}
Let $X_0$ be the total space of complex $\mathcal O(-1)\otimes
\mathbb C^n\rightarrow \mathcal O\otimes \mathbb C^{2n+1}\rightarrow
\mathcal O(1)\otimes \mathbb C^n$, $G_0$ be group
GL$_n\times$GL$_{2n+1}\times$GL$_n$/ $\mathbb C^\times$, $\rho_+$ be
the character $(1,0,-1)+\epsilon(n,-2n-1,0)$ for $\epsilon$ small
enough. Then $X_0$//$_{\rho_+}G_0$ is smooth. \label{smoothness of
X/0+ G}
\end{lemma}
\begin{proof}
For a stable complex \textbf K with respect ot $\rho_+$, we may
restricted it to the elliptic curve $E$, since Hom$(\textbf
K_E,\textbf K_E)$ is $\mathbb C$, the hypercohomology of \textbf
H$^2( \mathcal Hom^{\bullet}(\textbf K|_E,\textbf K|_E))$ is the
same as Ext$^2(\textbf K,\textbf K)$. Since $\textbf K|_E$ is exact
at the first term and the homological sheaf at the middle is a line
bundle with non-positive degree, it is quasi-isomorphic to $\mathcal
Q\rightarrow \overline{\mathcal L}^{\oplus n}$, where $\mathcal Q$
is locally free and $\mu_+(\mathcal Q)\leq 3 =\mu(\overline{\mathcal
L})$. Hence \textbf H$^2( \mathcal Hom^{\bullet} (\textbf
K|_E,\textbf K|_E))$ $=$ $0$. By a similar argument as that in
Corollary \ref{smoothness and irred property}, $X$//$_{\rho_+}G$ is
smooth.
\end{proof}
 By Proposition \ref{recollections from VGIT}, property 5,
since no curve is contracted, we have a birational map $T_w$:
$X_0$//$_{\rho_-}G_0$ $\dashrightarrow$ $X_0$//$_{\rho_+}G_0$, where
$X$//$_{\rho_-}G$ is Hilb$^nS$. As both varieties are smooth and
$T_w$ doesn't have exceptional locus, this is an isomorphism. Under
this isomorphism, the line bundle complex remains the same (since
they are stable on both sides). Moreover, due to the uniqueness of
the $S$-equivariant class, the $T_w$ image of an ideal complex
$\mathcal I_Z$ with $Z$ to be $n$ general distinct points
$p_1,\dots,p_n$ (by the term `general', we mean $\sigma^3(p_i)\neq
p_j$, $p_i\neq p_j$ for any $1\leq i,j\leq n$) is shown below.
\\\\
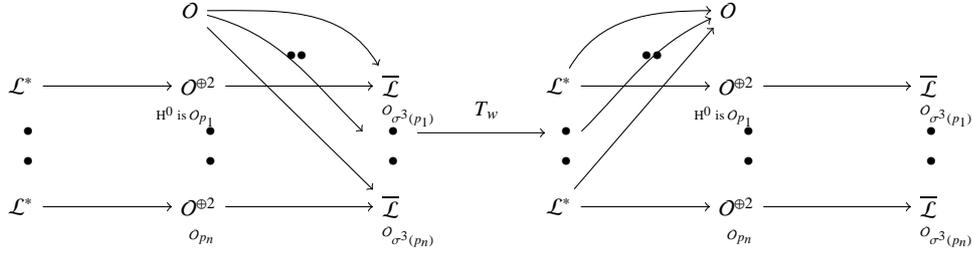
\begin{tikzpicture}[
%
  x=2.4cm,y=.2cm,
  font=\footnotesize,
  main/.style={draw,fill=yellow,inner sep=.5em},
  PP/.style={draw,fill=purple!40!blue!30,inner sep=1em},
  R/.style={draw,fill=purple!40!blue!30,inner sep=.5em},
  M/.style={draw,fill=green!80!yellow,inner sep=.5em},
  S/.style={anchor=east},
  V/.style={anchor=west},
  P/.style={anchor=center},
  F/.style={anchor=west}
  ]
\node[S] (C){}; \node[S] at ($(C)+(-3.1,9)$) (AL){$\mathcal O$};
\node[S] (C){}; \node[S] at ($(C)+(-3,4)$) (E1){$\mathcal
O^{\oplus2}$}; \node[S] at ($(C)+(-4,4)$) (E0){$\mathcal L^*$};
\node[S] (C){}; \node[S] at ($(C)+(-2,4)$) (E2){$\overline{\mathcal
L}$}; \node[S] (C){}; \node[S] at ($(C)+(-3,1)$) (){$\bullet$};
\node[S] at ($(C)+(-4,1)$) (){$\bullet$}; \node[S] (C){}; \node[S]
at ($(C)+(-2,1)$) (){$\bullet$}; \node[S] (C){}; \node[S] at
($(C)+(-3,-1)$) (){$\bullet$}; \node[S] at ($(C)+(-4,-1)$)
(){$\bullet$}; \node[S] (C){}; \node[S] at ($(C)+(-2,-1)$)
(){$\bullet$}; \node[S] (C){}; \node[S] at ($(C)+(-3,-4)$)
(G1){$\mathcal O^{\oplus2}$}; \node[S] at ($(C)+(-4,-4)$)
(G0){$\mathcal L^*$}; \node[S] (C){}; \node[S] at ($(C)+(-2,-4)$)
(G2){$\overline{\mathcal L}$}; \node[S] (C){}; \node[S] at
($(C)+(-3,2)$) (B1){\tiny{H$^0$ is $\mathcal O_{p_1}$}}; \node[S]
(C){}; \node[S] at ($(C)+(-1.8,2)$) (B2){\tiny{$\mathcal
O_{\sigma^3(p_1)}$}}; \node[S] (C){}; \node[S] at ($(C)+(-3,-6)$)
(D1){\tiny{$\mathcal O_{p_n}$}}; \node[S] (C){}; \node[S] at
($(C)+(-1.8,-6)$) {\tiny{$\mathcal O_{\sigma^3(p_n)}$}};

\node[S] at ($(C)+(-2.5,6)$) {$\bullet \bullet$};

\node[S] at ($(C)+(3,0)$) (C) {}; \node[S] at ($(C)+(-3.1,9)$)
(AR){$\mathcal O$};
 \node[S] at ($(C)+(-3,4)$) (ER1){$\mathcal O^{\oplus2}$};
\node[S] at ($(C)+(-4,4)$) (ER0){$\mathcal L^*$};
 \node[S] at ($(C)+(-2,4)$) (ER2){$\overline{\mathcal L}$};
 \node[S] at ($(C)+(-3,1)$) (){$\bullet$}; \node[S] at ($(C)+(-4,1)$) (){$\bullet$};
 \node[S] at ($(C)+(-2,1)$) (){$\bullet$};
 \node[S] at ($(C)+(-3,-1)$) (){$\bullet$}; \node[S] at ($(C)+(-4,-1)$) (){$\bullet$};
 \node[S] at ($(C)+(-2,-1)$) (){$\bullet$};
 \node[S] at ($(C)+(-3,-4)$) (GR1){$\mathcal O^{\oplus2}$}; \node[S] at ($(C)+(-4,-4)$) (GR0){$\mathcal L^*$};
 \node[S] at ($(C)+(-2,-4)$) (GR2){$\overline{\mathcal L}$};
 \node[S] at ($(C)+(-3,2)$) (B1){\tiny{H$^0$ is $\mathcal O_{p_1}$}};
 \node[S] at ($(C)+(-1.8,2)$) (B2){\tiny{$\mathcal O_{\sigma^3(p_1)}$}};
 \node[S] at ($(C)+(-3,-6)$) (D1){\tiny{$\mathcal O_{p_n}$}};
 \node[S] at ($(C)+(-1.8,-6)$) {\tiny{$\mathcal O_{\sigma^3(p_n)}$}};

\node[S] at ($(C)+(-3.5,6)$) {$\bullet \bullet$};

\node[S] at ($(-1.5,2.4)$) {$T_w$};

\draw[->] (E0) -- (E1); \draw[->] (E1) -- (E2); \draw[->] (G0) --
(G1); \draw[->] (G1) -- (G2); \draw[->] (AL) to [out =0,in =120]
(E2); \draw[->] (AL) -- (G2); \draw[->] (AL) to [out=-15,in=135]
(-2.3,1);

\draw[->] (ER0) -- (ER1); \draw[->] (ER1) -- (ER2); \draw[->] (GR0)
-- (GR1); \draw[->] (GR1) -- (GR2); \draw[->] (GR0) -- (AR);
\draw[->] (ER0) to [out =60, in =180] (AR);

\draw[->] (-1.1,1) to [out=45,in=-155] (AR);

\draw[->] (-2,0.9) -- (-1.3,0.9);


\end{tikzpicture}
\\\\
By writing a complex \textbf K in $X_0^{s,\rho_-}$ as $\mathcal
O(-1)\otimes H_{-1}\xrightarrow{I} \mathcal O\otimes
H_{0}\xrightarrow{J} \mathcal O(1)\otimes H_1$ with $I$ $=$ $xI_1$
$+$ $yI_2$ $+$ $zI_3$, $J$ $=$ $xJ_1$ $+$ $yJ_2$ $+$ $zJ_3$. Another
morphism $\tilde{T}_t$ from $X_0^{s,\rho_-}$  to $X_0^{s,\rho_+}$ is
defined as:
\begin{center}
$(I,J)$
$=$ ($xI_1$ $+$ $yI_2$ $+$ $zI_3$, $xJ_1$ $+$ $yJ_2$ $+$ $zJ_3$)
$\mapsto$ ($xJ^T_2$ $+$ $yJ^T_1$ $+$ $zJ^T_3$,$xI^T_2$ $+$ $yI^T_1$
$+$ $zI^T_3$).
\end{center}
\begin{lemma}
$\tilde{T}_t$ is well-defined and compatible with the $G_0$-action.
In addition, it extends to other quiver regions as
$\widetilde{T}_{t,k}:X^{s,\rho_-}_k\rightarrow X^{s,\rho_+}_{-k}$.
\label{lemma: tilde T is well defined}
\end{lemma}
\begin{proof}
1. Since $x$, $y$, $z$ satisfies the relations (\ref{relation
of x,y,z}), the image is really a complex. \\

2. The stability property is due to the duality.
$\tilde{T}_t(\textbf K)$ is a complex $\mathcal O(-1)\otimes
H^*_{1}\xrightarrow{I} \mathcal O\otimes H^*_{0}\xrightarrow{J}
\mathcal O(1)\otimes H^*_{-1}$. A subcomplex in $\tilde{T}_t(\textbf
K)$ is determined by  subspaces $(H'_1,H'_0,H'_{-1})$ in
$(H^*_1,H^*_0,H^*_{-1})$ those are compatible with
$\tilde{T}_t(I,J)$. Then $(H'^\perp_{-1},H'^\perp_0,H'^\perp_{1})$
in $(H_{-1},H_0,H_1)$ are compatible with $I$ and $J$, hence they
determine a subcomplex of \textbf K. Since
$\rho_+\cdot(h'_1,h'_0,h'_{-1})>0$ if and only if
$\rho_-\cdot(n_{-1}-h'_{-1},n_0-h_0,n_1-h_1)>0$,
$\tilde{T}_t(\textbf K)$ is $\rho_+$ stable.

\end{proof}
As $\tilde{T}_t$ maps a $G_0$-orbit to a $G_0$-orbit, it induces a
map from $X_0$//$_{\rho_-}G_0$ to $X_0$//$_{\rho_+}G_0$. We denote
this isomorphism between $X^{s,\rho_-}$ to $X^{s,\rho_+}$ by $T_t$.
This sets up the symmetry wall
crossing picture between the first and second quadrant.\\

Denote $T$ $:=$ $T_t\circ T_w$ by the automorphism of
$X_0$//$_{\rho_-}G_0$ $\simeq$ Hilb$^nS$. By the definition of
$T_t$, we have $T\circ T =$Id. The following statement shows that
when $n\geq 3$, the induced $T$-action on NS$_{\mathbb
R}$(Hilb$^nS)$ is non-trivial, i.e. the destabilizing wall on the
first quadrant destabilizes different points as those on the second
quadrant.
\begin{remark}
This involution $T$ is related to the Galois representation of the
symplectic resolution.
\end{remark}
\begin{proposition}
When $n\geq 3$, the automorphism $T$ induces a non-trivial action on
H$^2($ Hilb$^nS, \mathbb Z)$.   \label{left and right are different}
\end{proposition}
\begin{proof} When $n=3$, since the $\mathcal
O(-1)$-wall (respectively, $\mathcal O(1)[1]$-wall) is the first
wall on the left (right) of $t$-axis, it is enough to show that
these two walls destabilize different points on
  $X_0$//$_{\rho_+}G_0$. We study when an ideal sheave $\mathcal I_Z$ that can be written as the
kernel of $\mathcal O\rightarrow \oplus \mathcal O_{p_i}$ for $3$
general distinct points $p_1,p_2,p_3$ on $E$ is destabilized on the
$\mathcal O(-1)$-wall. Let the complex of $\mathcal I_Z[1]$ be
$\mathcal (L^*)^{\oplus 3}\rightarrow \mathcal O_E^{\oplus
7}\rightarrow \overline{\mathcal L}^{\oplus 3}$ as the cartoon on
the left. Write $\mathcal E$ for the kernel of $\mathcal O_E^{\oplus
7}\rightarrow \overline{\mathcal L}^{\oplus 3}$. As in the cartoon,
$\mathcal O_E^{\oplus 7}\rightarrow \overline{\mathcal L}^{\oplus
3}$ has four parts: $\mathcal O_E$ and three pieces of $\mathcal
O_E^{\oplus 2}\rightarrow \overline{\mathcal L}$. Each $\mathcal
O_E^{\oplus 2}\rightarrow \overline{\mathcal L}$ has kernel
$\overline{\mathcal L}^{-1}(\sigma^3(p_i))$ and cokernel $\mathcal
O_{\sigma^3(p_i)}$. The map from $\mathcal E$ to the direct sum of
the three pieces $\mathcal O_E^2\rightarrow \overline{\mathcal L}$,
has kernel $\mathcal O(-\sigma^3(p_1)-\sigma^3(p_2)-\sigma^3(p_3))$.\\

Since Hom$(\mathcal O(-1),\mathcal O(i)\otimes \mathbb C^{n_i})$'s
have dimensions 3,21,18, for $i=-1,0,1$ respectively,
\begin{align*}
&\text{Hom}(\mathcal O(-1), \mathcal I_Z)\neq 0
\\
\Leftrightarrow &\text{ the map from Hom}(\mathcal O(-1),\mathcal
O\otimes \mathbb C^{7})\text{ to Hom}(\mathcal O(-1),\mathcal
O(1)\otimes \mathbb C^{3})\text{ is not
surjective}\\
\Leftrightarrow &\text{ the map from Hom}(\mathcal L^*,\mathcal
O_E\otimes \mathbb C^{7})\text{ to Hom}(\mathcal
L^*,\overline{\mathcal L}\otimes
\mathbb C^{3})\text{ is not surjective}\\
\Leftrightarrow &\text{ Ext}^1(\mathcal L^*, \mathcal E)\neq 0\\
\Leftrightarrow & \mathcal
O(-\sigma^3(p_1)-\sigma^3(p_2)-\sigma^3(p_3)) \simeq \mathcal L^*.
\end{align*}
The last `$\Leftrightarrow$' is due to the short exact sequence
$0\rightarrow$ $\mathcal
O\big(-\sigma^3(p_1)-\sigma^3(p_2)-\sigma^3(p_3)\big)\rightarrow
\mathcal E\rightarrow$ $\oplus \overline{\mathcal
L}^{-1}(\sigma^3(p_i))$ $\rightarrow 0$. A similar argument shows
that $T_w(\mathcal I_Z[1])$, whose cartoon is on the right of the
previous picture, has non-zero morphism to $\mathcal O(1)[1]$ if and
only if $\mathcal O(p_1+p_2+p_3) \simeq \overline{\mathcal L}$.
Hence $\mathcal O(-1)$ has non-zero morphism to $T(\mathcal I_Z)$ if
and only if $\mathcal O(p_1+p_2+p_3) \simeq \overline{\mathcal L}$.
Since $\overline{\mathcal L}(p_i)=\mathcal L^* (\sigma^3(p_i))$ and
$\sigma$ has infinite order, the locus that is contracted by the
$\mathcal O(-1)$-wall and that is contracted by the $\mathcal O(1)[1]$-wall are different.\\

When $n\geq 4$, we do the induction on $n$. Assume the $n-1$ case is
done, then a line bundle $\mathcal I$ with $(r,c_1,\chi)$ $=$
$(1,0,1-(n-1))$ is destabilized by $\mathcal O(-1)$, and $T(\mathcal
I)$ is not destabilized by $\mathcal O(-1)$. Consider the morphism
$\mathcal O(-1)\rightarrow \mathcal I$ restricted on $E$, the
cokernel is a torsion sheaf of length $3$, let $\mathcal O_p$ be a
quotient of the torsion sheaf. Then $\mathcal O(-1)$ has a non-zero
map to the kernel $\mathcal I'$ of $\mathcal I\rightarrow$ $\mathcal
O_p$. Yet $T(\mathcal I')$ is the kernel of $T(\mathcal
I)\rightarrow \mathcal O_q$ for some $q\in
E$, Hom$(\mathcal O(-1),T(\mathcal I'))$ $=$ $0$.\\

 Since for any destabilize sequence
$\mathcal O(-1)\rightarrow \mathcal I'\rightarrow \mathcal I''$. The
extension sheaf by $\mathcal O(-1)$ and $\mathcal I''$ is a vector
bundle if and only if for any non-zero numbers $(l_1,l_2,l_3)$ $\in$
$\mathbb C^3$ on $E$, $aI_x$ $+$ $bI_y$ $+$ $cI_z$ is injective i.e
$l_1I^T_x$ $+$ $l_2I^T_y$ $+$ $l_3I^T_z$ is surjective. For generic
choice of Hom$(\mathcal O(-2)\otimes \mathbb C^{n_{-1}},\mathcal
O(-1))$ and $\mathcal I''$, $l_1I^T_x$ $+$ $l_2I^T_y$ $+$ $l_3I^T_z$
is injective since for generic $\mathcal I''$ the cokernel of
$xI''^T_x$ $+$ $yI''^T_y$ $+$ $zI''^T_z$ restricts on $E$ is the
direct sum of some skyscraper sheaves of distinct points. Hence on
the locus that are destabilized by $\mathcal O(-1)$, the set of
vector bundles is dense. Therefore there exists a vector bundle that
is destabilized by $\mathcal O(-1)$ while $T(-)$ of it is not
destabilized by $\mathcal O(-1)$. The induction accomplishes.
\end{proof}

Combining Theorem \ref{left half upper plane's main theorem in the
body} and Proposition \ref{left and right are different}, we get our
main result.
\begin{theorem}
When $n\geq 3$, the positivity cone of Hilb$^nS$ is symmetric. Each
side stable base locus decomposition walls are one to one
corresponding to the semicirclar actual walls on the first and
second quadrant of Bridgeland stability conditions. \hfill
$\square$\label{main theorem in the body}
\end{theorem}
\section{Examples}
\begin{example}
Given $n$, when $\sqrt{2n}>k\geq 0$, in the quiver region of
$\mathcal A(-k)$, the character $\overrightarrow{\rho}$ is given by:
\begin{center}$\frac{ts}{\frac{t^2}{2}+2n-s^2}((s+k+1)^2-\frac{t^2}{2},-(s+k)^2+\frac{t^2}{2},(s+k-1)^2-\frac{t^2}{2})+t(s+k+1,-s-k,s+k-1).
$\label{character rho of s,t}
\end{center}
\end{example}
When $t$ tends to $0$, the character $\overrightarrow{\rho}_{s,t,k}$
of $G/\mathbb C^\times$ is up to a scalar given by:
\[(s^2(k+1)+s(2n+(k+1)^2)+2n(k+1),-s^2k-s(2n+k^2)-2nk,(s^2(k-1)+s(2n+(k-1)^2)+2n(k-1)).\]
When $s$ decreases from $-k+1$ to $-k-1$, the character decreases
from $(1,-\frac{n_{-1}}{n_{0}},0)$ to  $(0,\frac{n_{1}}{n_{0}},-1)$.
In particular, when $s$ is $-k$ and $t$ tends to $0$, up to a scalar
$\overrightarrow{\rho}_{-k,0+,k}$ is $(n_1,0,-n_{-1})$, it
corresponds to the destabilizing walls with
type $(0,1,0)$, as a sheaf it is just $\mathcal O(-k)[1]$.\\

Given an integer $-\sqrt{2n}$ $<$ $k$ $\leq$ $0$, for $-k-1<s<-k+1$,
let $A_k$ and $B_k$ be the line bundles (divisors) on
$X$//$_{\overrightarrow{\rho}_{s,0+,k}}G$ that compose with the
$G$-principal bundle with characters $(1,*,0)$ and $(0,*,-1)$
respectively. Then when $s$ is between two integers $-k-1$ and $-k$,
there are four divisors $A_k$, $B_k$, $A_{k+1}$ and $B_{k+1}$. When
quiver region only contains flip-type bi-rational morphism, by the
Remark \ref{remark on the divisor glueing in flip case}, these
divisors satisfy the relation:
\[
c_k
\begin{bmatrix}
A_{k+1}\\
B_{k+1}
\end{bmatrix}
=\begin{bmatrix} 2n-k(k+1) & 2n-k(k-1)\\
-2n+(k+1)(k+2) & 3(2n-(k-1)(k+2))
\end{bmatrix}
\begin{bmatrix}
A_k\\B_k\end{bmatrix}, \tag{$\triangle$} \label{glueing divisors
formula}
\]
where $c_k$ is a constant only depend on $k$. Furthermore $A_k$
$\sim$ $B_{k-2}$, where $\sim$ means equal up to a scalar.

\begin{proposition}
Let the notations $A_k$, $B_k$ be as above.  Assume $A_1$ $\sim$
$H$, $B_0$ $\sim$ $A_2$ $\sim$ $(n-1)H-\frac{\Delta}{2}$, $A_3$
$\sim$ $\frac{n-1}{2} H-\frac{\Delta}{2}$, then the divisor at
$(s,0+)$ is $(-\frac{2n+s^2}{2s}-\frac{3}{2})H-\frac{\Delta}{2}$ up
to scalar. In an other word, the destabilizing semicircle wall on
the Bridgeland stability condition space with center
$-m-\frac{3}{2}$ corresponds to the divisor $mH-\frac{\Delta}{2}$.
\label{proposition: the explicit formula of 1-1 walls}
\end{proposition}
\begin{proof}
First of all, we show that $A_k$ and $B_k$ are
$(2n+(k-1)(k-4))H-(k-1)\Delta$ and $(2n+(k-2)(k+1))H-(k+1)\Delta$
respectively up to a same scalar.\\
When $k=1$, we may assume that $A_1 = 2nH$, $B_1=
b_1((n-1)H-\Delta)$, $A_2 =a_2((n-1)H-\frac{\Delta}{2})$. By the
equation (\ref{glueing divisors formula}), we have
\[ A_1(2n-2)+2nb_1(\frac{n-1}{2}H-\Delta)\sim
(n-1)H-\frac{\Delta}{2}.\] This implies $b_1 =2$. By the equation
(\ref{glueing divisors formula}) and induction on $k$, we get $A_k$
and $B_k$.\\

At a point $(s,0+)$, the character $\rho_{s,0+,k}$ is given in
Example \ref{character rho of s,t}. As
\[\rho_{s,0+,k} =
-f(n,s,k-1)(0,\frac{n_1}{n_0},-1)+f(n,s,k+1)(1,-\frac{n_{-1}}{n_0},0),\]
where $f(n,s,k)$ $=$ $k(2n+s^2)+s(2n+k^2)$. The divisor at $(s,0+)$
is up to a scalar given by:
\begin{align*}
& -f(n,s,k-1)B_k+f(n,s,k+1)A_k\\
\sim &
-f(n,s,k-1)((2n+(k-2)(k+1))H-(k+1)\Delta)+f(n,s,k+1)((2n+(k-4)(k-1))H-(k-1)\Delta)\\
= & 2(2n-(k-1)(k+1))(2n+s^2+3s)H +2s(2n-(k-1)(k+1))\Delta\\
= & -2s(2n-(k-1)(k+1))\big(
(-\frac{2n+s^2}{2s}-\frac{3}{2})H-\frac{1}{2}\Delta\big).
\end{align*}
\end{proof}

\subsection{Destabilizing Walls}
To compute the ratio of each stable decomposition wall on the
Neron-Severi space, we only need compute all the ratio
$\rho_{s,0+,k}$'s on the destabilizing chamber walls. We may look at
each $\mathcal A(k)$ quiver region to search candidates type of
subcomplex that may destabilize a stable complex \textbf K with type
$(n-\frac{(k-1)k}{2},2n-k^2+1,n-\frac{k(k+1)}{2})$. \\

For each quiver region, we only need consider the wall whose right
bound is in $(-k-1,-k)$. Suppose the character $\overrightarrow
{\rho}$ gives an actual wall, then there is a destabilizing
sequence: \textbf K$''\rightarrow$ \textbf K $\rightarrow $ \textbf
K$'$ with \textbf K$'$ stable. Let the type of \textbf K$''$ be
$(a+l,2a+r+l,a)$, then \textbf K$'$ has type $(A,A+C-s,C)$ $=$
$(n-\frac{(k-1)k}{2},2n-k^2+1,n-\frac{k(k+1)}{2})$ -
$(a+l,2a+r+l,a)$. To achieve an efficient logarithm, we need some
restrictions on the candidate type $(a+l,2a+r+l,a)$.
\begin{lemma}
Let $a$, $r$, $l$ be as discussed before, then they satisfy the
following inequalities:
\begin{align}
& a+l \geq 0;\label{inequality: a+l geq 0}\\
& (A-C)^2-s(A+C-s)-2s^2+1  \geq  0;\label{chi2 vanishing inequality}\\
& \frac{k+1}{n-\frac{k(k+1)}{2}}a-r < l < \frac{k}{n-\frac{k(k+1)}{2}}a ;\label{ratio falls into inequality}\\
& l+r\leq a\text{; if }r\geq 2\text{, then }2a \geq 3(r+l)\text{ or
the type is}
(0,3,1); \dots \label{031 type inequality}\\
& k  >  \sqrt{(r-1)(2n+r-1)/r}-1. \label{inequality: boundary for r}
\end{align}
\end{lemma}
\begin{proof}
Inequality (\ref{chi2 vanishing inequality}) is a consequence of
Lemma \ref{lemma: ext2 vanishing for stable factors}. Since
Ext$^2(\textbf K',\textbf K')$ $=$ $0$ and $\textbf K'$ is stable,
we have
\[ \chi(\textbf K',\textbf K')\leq \text{dimHom}( \textbf K',\textbf K') = 1.\]
On the other hand,
\begin{align*}
\chi(\textbf K',\textbf K') &=
\text{dimHom}^0(\textbf K',\textbf K')-\text{ dimHom}^1(\textbf
K',\textbf K')+\text{ dimHom}^2(\textbf K',\textbf K')\\
&= (A^2+(A+C-s)^2 +C^2)-(3A(A+C-s)+3C(A+C-s)) +6AC\\
&= -(A-C)^2+s(A+C-s)+2s^2.
\end{align*}
Inequality (\ref{ratio falls into inequality}): by formula
\ref{character rho of s,t}, the boundary $-k-1$ and $-k$ corresponds
to characters $\overrightarrow{\rho}_{-k-1,0+,k}$ $\sim$
$(0,n-\frac{(k+1)k}{2},-(2n-k^2+1))$ and
$\overrightarrow{\rho}_{-k,0+,k}$ $\sim$
$(n-\frac{(k+1)k}{2},0,-(n-\frac{k(k-1)}{2}))$. We have
\begin{center}
$ (a+l,2a+r+l,a)\cdot \overrightarrow{\rho}_{-k-1,0+,k}$ $>$ $0$;\\
$ (a+l,2a+r+l,a)\cdot \overrightarrow{\rho}_{-k,0+,k}$ $<$ $0$.
\end{center}
Plug in the values, we get the two boundaries for $l$.\\

Inequality (\ref{031 type inequality}): if $2a+r+l>3a$, we may
consider the intersection of ker$J''_x$,  ker$J''_y$ and ker$J''_z$,
then \textbf K$''$ contains $(0,1,0)$ type sub complex, $\textbf K$
is already destabilized at a previous wall.\\

The formula \ref{inequality: boundary for r} is implied by \ref{chi2
vanishing inequality} and \ref{ratio falls into inequality}. Write
the inequality \ref{chi2 vanishing inequality} in terms of $a$, $l$
and $r$: one has
\[L:=(k-l)^2-(r-1)(2n-k^2+1-2a-l-r)-2(r-1)^2+1 \geq 0.\]
When $r\leq 1$, the inequality holds obviously. We may assume $r\geq
2$. When $l\geq k+2-r$, we have $L$ $\leq$ $(r-2)^2$ $-$ $2(r-1)^2$
+ $1$ $<$ $0$, hence $l\leq k+1-r$.\\

By the first part of \ref{ratio falls into inequality}, we have:
\[(k-l)^2-(r-1)(2n-k^2+1-2\frac{n-\frac{k(k+1)}{2}}{k+1}(l+r)-l-r)-2(r-1)^2+1 > 0.\]
Now the left side is a binomial of $l$ with leading coefficient $1$.
If an $l$ $\in$ $(1-r , k+1-r)$ satisfies the inequality, then
either $1-r$ or $k+1-r$  satisfies it. Plug in $l$ $=$ $k+1-r$, the
inequality always fails. Hence it holds for $l = 1-r$.
\begin{align*}
&(k+r-1)^2-(r-1)(2n-k^2-2\frac{n-\frac{k(k+1)}{2}}{k+1})-2(r-1)^2+1\geq
0\\
\Leftrightarrow & rk^2+(r-1)k+1\geq (r-1)(r-1+\frac{k}{k+1}2n)\\
\Rightarrow & \frac{k+1}{k}rk(k+1)\geq (r-1)(r-1+2n)\\
\Rightarrow & k  >  \sqrt{(r-1)(2n+r-1)/r}-1.
\end{align*}

\end{proof}

\bibliographystyle{abbrv}
\bibliography{MMPref}

\begin{thebibliography}{10}

\bibitem{ABCH}
D.~Arcara, A.~Bertram, I.~Coskun, and J.~Huizenga.
\newblock The minimal model program for the {H}ilbert scheme of points on
  {$\Bbb{P}^2$} and {B}ridgeland stability.
\newblock {\em Adv. Math.}, 235:580--626, 2013.

\bibitem{ATV1}
M.~Artin, J.~Tate, and M.~Van~den Bergh.
\newblock Some algebras associated to automorphisms of elliptic curves.
\newblock In {\em The {G}rothendieck {F}estschrift, {V}ol.\ {I}}, volume~86 of
  {\em Progr. Math.}, pages 33--85. Birkh\"auser Boston, Boston, MA, 1990.

\bibitem{ATV2}
M.~Artin, J.~Tate, and M.~Van~den Bergh.
\newblock Modules over regular algebras of dimension {$3$}.
\newblock {\em Invent. Math.}, 106(2):335--388, 1991.

\bibitem{AV}
M.~Artin and M.~Van~den Bergh.
\newblock Twisted homogeneous coordinate rings.
\newblock {\em J. Algebra}, 133(2):249--271, 1990.

\bibitem{Br}
T.~Bridgeland.
\newblock Stability conditions on triangulated categories.
\newblock {\em Ann. of Math. (2)}, 166(2):317--345, 2007.

\bibitem{Br2}
T.~Bridgeland.
\newblock Stability conditions on {$K3$} surfaces.
\newblock {\em Duke Math. J.}, 141(2):241--291, 2008.

\bibitem{CH}
I.~Coskun and J.~Huizenga.
\newblock Interpolation, bridgeland stability and monomial schemes in the
  plane.
\newblock {\em preprint,
  http://homepages.math.uic.edu/~coskun/monomialBaseLocus.pdf}.

\bibitem{SGA}
P.~Deligne.
\newblock {\em Cohomologie \'etale}.
\newblock Lecture Notes in Mathematics, Vol. 569. Springer-Verlag, Berlin,
  1977.
\newblock S{\'e}minaire de G{\'e}om{\'e}trie Alg{\'e}brique du Bois-Marie SGA
  4${1{\o}er 2}$, Avec la collaboration de J. F. Boutot, A. Grothendieck, L.
  Illusie et J. L. Verdier.

\bibitem{DoHu}
I.~V. Dolgachev and Y.~Hu.
\newblock Variation of geometric invariant theory quotients.
\newblock {\em Inst. Hautes \'Etudes Sci. Publ. Math.}, (87):5--56, 1998.
\newblock With an appendix by Nicolas Ressayre.

\bibitem{Gi}
V.~Ginzburg.
\newblock Lectures on nakajima's quiver varieties.
\newblock {\em preprint}.

\bibitem{Gr}
A.~Grothendieck.
\newblock El\'ements de g\'eom\'etrie alg\'ebrique, chapters iii and iv.
\newblock {\em Inst. Hautes \'Etudes Sci. Publ. Math.}, 11, 1961.

\bibitem{Hi}
N.~Hitchin.
\newblock Deformations of holomorphic {P}oisson manifolds.
\newblock {\em Mosc. Math. J.}, 12(3):567--591, 669, 2012.

\bibitem{Ki}
A.~D. King.
\newblock Moduli of representations of finite-dimensional algebras.
\newblock {\em Quart. J. Math. Oxford Ser. (2)}, 45(180):515--530, 1994.

\bibitem{Ma}
E.~Macr{\`{\i}}.
\newblock Stability conditions on curves.
\newblock {\em Math. Res. Lett.}, 14(4):657--672, 2007.

\bibitem{Nak}
H.~Nakajima.
\newblock Heisenberg algebra and {H}ilbert schemes of points on projective
  surfaces.
\newblock {\em Ann. of Math. (2)}, 145(2):379--388, 1997.

\bibitem{NS}
T.~A. Nevins and J.~T. Stafford.
\newblock Sklyanin algebras and {H}ilbert schemes of points.
\newblock {\em Adv. Math.}, 210(2):405--478, 2007.

\bibitem{SV}
J.~T. Stafford and M.~van~den Bergh.
\newblock Noncommutative curves and noncommutative surfaces.
\newblock {\em Bull. Amer. Math. Soc. (N.S.)}, 38(2):171--216, 2001.

\bibitem{Th}
M.~Thaddeus.
\newblock Geometric invariant theory and flips.
\newblock {\em J. Amer. Math. Soc.}, 9(3):691--723, 1996.

\end{thebibliography}
\end{document}